\RequirePackage{amsmath}
\documentclass[onecolumn]{svjour3} 
\smartqed  
\usepackage{graphicx}
\usepackage{subfig}
\usepackage[utf8]{inputenc}
\usepackage{dsfont}
\usepackage{relsize}
\usepackage{mathptmx}
\usepackage{url}
\usepackage{amsmath}
\usepackage{amssymb}
\usepackage{amsfonts}

\newcommand{\R}{\mathbb{R}}

\usepackage[usenames, dvipsnames]{color}

\title{Robust topology optimization using a posteriori error estimator for the finite element method}
\author{Pimanov, Vladislav \and Oseledets, Ivan }
\institute{V. Pimanov \and I. Oseledets \at Skolkovo Institute of Science and Technology, Nobel St. 3, Moscow, Russia \and
I. Oseledets \at Institute of Numerical Mathematics of Russian Academy of Sciences, Gubkin St. 8, Moscow, Russia.}
\journalname{Struct Multidisc Optim}
\date{}
\begin{document}
\maketitle
\begin{abstract}
In our work, we consider the classical density-based approach to the topology optimization. We propose to modify the discretized cost functional using a posteriori error estimator for the finite element method. It can be regarded as a new technique to prevent checkerboards. It also provides higher regularity of solutions and robustness of results.
\keywords{topology optimization \and heat conduction \and finite element method
\and a posteriori error estimation \and checkerboards.}
\end{abstract}

\begin{section}{Introduction}\label{sec:introduction}
In our study, we consider the classical \textit{density-based} approach to topology optimization problems which consists in distributing of material inside a fixed domain and assumes the material is modeled as a piece-wise constant on a fixed finite element mesh function.
The topology optimization problem is reduced to the minimization of the \textit{cost functional} on the specified set of \textit{admissible designs}.
Computation of the cost functional for any fixed design requires a boundary value problem (BVP) to be solved, and we actually deal with its numerical approximation. 
When the finite element method (FEM) is used, it can be shown (Section \ref{section:relation}) that for many problems encountered in practice including the model problem considered in our work, the
true value of the cost functional is always greater than its discrete value,
and the error of the functional is straightforwardly determined by the error of the finite element solution of the underlying BVP.
The main difficulty is that discontinuity and strong heterogeneity of
the coefficients, relevant to topology optimization problems, often lead to a poor approximation of solutions in the standard finite element subspaces of piece-wise polynomial functions.
So, small values of the discrete cost functional often do not
lead to small values of the true cost functional.
In particular, checkerboard-like designs \cite{sigmund-numinstabilities-1998,diaz-checkerboard-1995}
clearly demonstrate such a ``false minima'' problem, since their formation is exactly due to a poor numerical modeling by lower order finite elements
and can not be interpreted as a kind of optimal porous microstructure \cite{sigmund-topopt-1994,diaz-checkerboard-1995}.
We propose a new technique that builds upon the ideas of \cite{ovchinnikov-apostreg-2017}.
During the minimization process, we also take into account
the FEM error using a posteriori error estimator.
We modify the discrete cost functional by an additional correction term that penalizes designs with a large FEM error.
Specifically, it can be regarded as a new technique to prevent checkerboards. 
In a broader sense, it intends to avoid ``false minima'' and to provide robust results.

Main contributions of our paper are:
\begin{itemize}
 \item We show the relation between true and discrete cost functionals through the FEM approximation error by a new interpretation of classical results of the FEM theory (Section \ref{section:relation})
 \item In Section \ref{section:modification}, we propose the modification of the discrete cost functional, which is based on a posteriori error estimator presented in Section \ref{section:apost}
 \item For the heat conduction model problem, we demonstrate that the minimization of this modified cost functional prevents formation of checkerboards and provides robust results (Section \ref{section:numexp}).
\end{itemize}

\end{section}

\begin{section}{Heat conduction model problem}\label{section:modelproblem}

In our paper, the model problem and results are expounded
with respect to the two-dimensional heat conduction problem.
From the mathematical point of view,
it is similar to the problem of the compliance optimization of the variable thickness sheet \cite{rossow-vts-1973,bendose-topoptbook-2013}, which is very well studied in the field of structural design.
The choice of such a model problem is primarily justified by its simplicity, even so, it is sufficient to illustrate the basic concepts.

The optimization task is to find optimal distribution of isotropic material inside a given domain to get the design with the maximal thermal response.
The design variable is the coefficient of thermal conductivity, and
the cost functional that we consider is the thermal compliance.
Let $\Omega \subset R^2$ be the polygonal Lipschitz domain with
boundary $\partial \Omega = \overline{\Gamma_u \cup \Gamma_n},\; \Gamma_u \cap \Gamma_n = \emptyset$, where
zero temperature is prescribed along the boundary $\Gamma_u$, and zero heat flux is prescribed along the boundary $\Gamma_n$.
We define the \textit{solution space}, denoted as $\mathcal{H}$, which is a subspace of the usual Sobolev space $H^1(\Omega)$:
\begin{equation*}
 \mathcal{H} = \{u \in H^1(\Omega)\;\big|\; u = 0\; \text{on}\; \Gamma_u\}.
\end{equation*} 
We consider the following set of admissible designs, denoted as $K_{ad}$, that admits intermediate values of the coefficients:
\begin{equation*}\label{continuousdesigns}
  K_{ad}(\Omega) = \big\{k \in L^{\infty}(\Omega)\;\big|\; \gamma \leq k \leq 1\;a.e.\;in\;\Omega;\;\int_{\Omega}k = \mathbf{V}\big\},
\end{equation*}
where $\mathbf{V}$ is the volume constraint, and $0 < \gamma \ll 1$ represents the conductivity of an ersatz material.
The cost functional, denoted as $\Phi(k)$, is the functional of the temperature distribution $u = u(k) \in \mathcal{H}$, hence its computation for any fixed design $k \in K_{ad}$ requires solving the \textit{underlying} boundary value problem (BVP). The Topology Optimization problem looks as follows:
\begin{equation}\label{problem:topopt}\tag{\textbf{TO}}
\left \{
\begin{aligned}
 &\textrm{{\large Minimize}} \;\; \Phi(k),\\
 (k,&u) \in K_{ad}\times\mathcal{H}\\
  &\textrm{{\large Subject to:}}\\ &\Phi(k) = \ell(u(k)),\\
   &a_{k}(u,v)=\ell(v),\;\forall v \in \mathcal{H},
\end{aligned}
\right.
\end{equation}
where $a_k(\cdot,\cdot):\mathcal{H}\times \mathcal{H} \to \mathbb{R}$ is a coercive symmetric continuous bi-linear form, associated with the design $k \in K_{ad}$,
and $\ell(\cdot):\mathcal{H} \to \mathbb{R}$ is a bounded linear form. In the case of heat conduction, these forms are given as follows:
\begin{equation}\label{forms}
\begin{aligned}
 a_k(u,v) &= \int_{\Omega} k^p \nabla u \cdot \nabla v,\\
 l(v) &= \int_{\Omega} f v,
\end{aligned}
\end{equation}
where $f \in L_2(\Omega)$ is a heat source, and $p \geq 1$ is a penalization factor that penalizes intermediate values of the coefficients, following the classical SIMP (Solid Isotropic Material with Penalization) approach \cite{bendsoe-simporigin-1989,rozvany-simpabbr-1992}.

It is a known fact that the \eqref{problem:topopt} problem generally has no solution when $p > 1$ \cite{bendose-topoptbook-2013}.
General ideas to deal with non-existence of solution
are to reduce the set of admissible designs by some sort of global
or local restrictions on the variation of the coefficients \cite{sigmund-numinstabilities-1998}.
For example, perimeter constrained \cite{ambrosio-perimetercontrol-1993,haber-perimetercontrol-1996}
or slope constrained \cite{niordson-slopecontrol-1983} sets of admissible designs can be considered.
In practice, regardless of whether a problem possesses a well-posed continuum formulation, designs are always discretized, and the existence issue does not arise in the case of finite dimensionality.
In the classical approach, the optimization model and the finite element model are strongly coupled
in the sense that the designs are approximated by the functions which are
piece-wise constant on the same FEM mesh that is used for solving the underlying BVP.
There are several fundamental theoretical studies dedicated to this approach. For example,
the convergence results for the variable thickness sheet problem without penalization are presented
in \cite{petersson-vts-1999}, and the convergence study in the case of slope
constrained set of admissible designs is considered in \cite{petersson-slope-1998}. 

In our work, we distinguish the \textit{model} grid, denoted as $M^H$, consisting of ground elements intended for the designs representation, from the \textit{computational} grid, denoted as $T_h$, intended for the temperature field approximation. Thus, the index $H$ denotes the characteristic model size,
when the index $h$ denotes the actual FEM mesh size.
Such notation is motivated by the desire to be able to refine the FEM mesh for a fixed problem.
We denote the set of piece-wise constant on $M^H$ functions as $K^H$
and define the discretized set of admissible designs:
\begin{equation*}
 K_{ad}^H = K_{ad} \cap K^H.
\end{equation*}
Then the discretized \eqref{problem:topopt} problem looks as follows:
\begin{equation}\label{problem:topopt-discrete}\tag{$\mathbf{TO}^H$}
\left \{
\begin{aligned}
 &\textrm{{\large Minimize}} \;\; \Phi(k),\\
 (k,&u) \in K^H_{ad}\times\mathcal{H}\\
  &\textrm{{\large Subject to:}}\\ &\Phi(k) = \ell(u(k)),\\
   &a_{k}(u,v)=\ell(v),\;\forall v \in \mathcal{H}.
\end{aligned}
\right.
\end{equation}
This \eqref{problem:topopt-discrete} problem is the one we actually want to study.
The lack of well-posedness of the \eqref{problem:topopt} problem leads to the
\textit{mesh-dependency} phenomena \cite{sigmund-numinstabilities-1998},
when the \eqref{problem:topopt-discrete} problem has qualitatively different solutions for different model grids,
so we can not discuss any convergence when $H$ goes to zero.
However, for any fixed model grid $M^H$, it can be considered as a completely independent and rather complicated task.

When we solve the \eqref{problem:topopt-discrete} problem in practice, we always deal with a 
discrete approximation of the true cost functional:
\begin{equation*}\label{discretecostfunc}
 \Phi_h(k) = \ell(u_h(k)),
\end{equation*}
where $u_h(k) \in \mathcal{H}_h$ is a finite element approximation of the true solution $u(k) \in \mathcal{H}$,
and $\mathcal{H}_h \subset \mathcal{H}$ is a finite element subspace.
In our work, we consider customary Lagrange quadrilateral finite elements with standard conforming piece-wise bi-linear and piece-wise bi-quadratic approximations:
\begin{equation}\label{fem-space}
  \mathcal{H}^l_h(T_h) = \; \Big\{v_h \in C(\overline{\Omega}) \; \Big| \; v_h |_T \in Q_l(T), \; \forall T \in T_h\Big\}\cap\mathcal{H}, \;\; l = \{1,2\}.
\end{equation}
The finite element discretization of the \eqref{problem:topopt-discrete} problem
with the discretized solution space looks as follows:
\begin{equation}\label{problem:topopt-discrete-fem}\tag{$\mathbf{TO}^H_h$}
\left \{
\begin{aligned}
 &\textrm{{\large Minimize}}\;\; \Phi_h(k),\\
 (k,&u_h) \in K^H_{ad}\times\mathcal{H}_h\\
  &\textrm{{\large Subject to:}}\\
  &\Phi_h(k) = \ell(u_h(k)),\\
  &a_{k}(u_h,v_h)=\ell(v_h),\;\forall v_h \in \mathcal{H}_h.
\end{aligned}
\right.
\end{equation}
Such separation of the \eqref{problem:topopt-discrete-fem} and the \eqref{problem:topopt-discrete} problems is primarily motivated by the desire to clearly designate that, in fact, we are interested 
in solving the \eqref{problem:topopt-discrete} problem but not the \eqref{problem:topopt-discrete-fem} problem, as it can be misunderstood in the case of a coupled
discretization of the solution and the admissible designs spaces on the same mesh.
We show in Section \ref{subsec:appr_problem} that the
approximation properties of the standard finite element subspaces \eqref{fem-space} are often not satisfactory
to ensure the solution of the \eqref{problem:topopt-discrete-fem} problem
to be close to the solution of the \eqref{problem:topopt-discrete} problem.
In the following Section we present the corresponding result of our study.
\end{section}

\begin{section}{Relation with the FEM error}\label{section:relation}
 Since the bi-linear form $a_k(\cdot,\cdot)$ is symmetric and coercive, it defines the energy inner product $(\cdot,\cdot)_a$
 with the corresponding energy norm $\|\cdot\|_a = a_k(\cdot,\cdot)^{1/2}$.
 We have the following interpretation of the classical results
 from the FEM theory (the Corollary of Theorem 1.1 in \cite{strang-femtheory-1973}) that establishes the relation between the true and the discrete cost functionals:
 \begin{theorem}\label{maintheorem}
 
 For an arbitrary designs $k \in K_{ad}$, let $u(k) \in \mathcal{H}$ be the true solution of the underlying BVP and $u_h(k) \in \mathcal{H}_h$ be its finite element approximation, then we have:
  \begin{equation*}
   \ell(u(k)) = \ell(u_h(k)) + \|u(k) - u_h(k)\|_a^2,
  \end{equation*}
  or the same in the context of the considered topology optimization problem:
  \begin{equation*}
  \Phi(k) = \Phi_h(k) + \|u(k) - u_h(k)\|_a^2.
  \end{equation*}
 \end{theorem}
 The following apparent corollary holds:
 \begin{corollary}
    $\Phi(k) \geq \Phi_h(k)$.
 \end{corollary}
Thus, small values of $\Phi_h(k)$ can be achieved due to a large FEM error.
We refer this phenomena as the ``false minima'' problem.
In particular, checkerboards clearly demonstrate the ``false minima'' problem since their formation is explained by a poor numerical modeling and their optimality is artificial.

A natural way to get checkerboard-free designs is to improve approximation properties of the discrete solution space $\mathcal{H}_h$.
For example, using higher order finite elements in each ground element or more than one finite element per
ground element usually helps to avoid checkerboards \cite{rozvany-topoptbook-2014,sigmund-lastreview-2013}. 
Another approach is to use special-type finite elements.
In the field of topology optimization, non-conforming finite elements providing checkerboard-free
results were studied in \cite{jang-nonconforming-2003,jang-nonconforming-2005}.
\begin{remark}
Theorem 1 is formulated in the terms of general forms $a_k(\cdot,\cdot)$ and $\ell(\cdot)$, so it holds for a wide class of all self-adjoint (i.e. the cost functional is strongly connected with the right-hand side: $\Phi(k) = \ell(u(k))$) topology optimization problems, where design variables are the coefficients of linear elliptic equations.
 Moreover, it can be generalized to the case of non-self-adjoint problems.
 For an arbitrary bounded linear cost functional $\Phi(\cdot):\mathcal{H} \to \R$, it holds that $\Phi(k) = \Phi_h(k) + (u(k)-u_h(k),z(k))_a$, where $z(k) \in \mathcal{H}$ is the solution of the corresponding adjoint BVP:
 \begin{equation*}
 	a_k(z,v) = \Phi(v), \; \forall v \in \mathcal{H}.
 \end{equation*}
\end{remark}
\end{section}

\begin{section}{FEM convergence results and quasi-monotonicity condition}\label{section:convergence}
The aim of this Section is to show how the FEM performs in the case of piece-wise constant coefficients relevant to topology optimization problems.
All the results presented in this Section can be found in details in the exhaustive theoretical 
study \cite{petzoldt-phd-2001}.
Generally, the asymptotic convergence rate of the FEM depends on the global regularity of the true solution $u \in \mathcal{H}$ and on the approximation properties of the finite element subspace $\mathcal{H}_h$.
We discuss the regularity of solutions using Sobolev spaces of fractional order $H^s(\Omega)$, $s \in \R$ as defined in \cite{adams-sobolev-2003} and denote its seminorm as $|\cdot|_{H^s(\Omega)}$.

Regularity results from \cite{jochmann-regularity-1999} show that, for an arbitrary design  $k \in K_{ad}^H$, it holds that $u(k) \in H^{1+s}(\Omega)$ for a certain $s > 0$.
In the case of a uniform computational grid $T_h$,
the approximation error can be measured in the terms of the grid size $h$. 
 For the piece-wise bi-linear finite element solution $u_h(k) \in \mathcal{H}_{h}^1$, we have:
\begin{equation}\label{convergence}
   \|u - u_h\|_a^2 \leq Ch^{2s}|u|_{H^{1+s}(\Omega)}^2,
\end{equation}
where the constant $C$ only depends on shape regularity of $T_h$.
The bad news is that the regularity parameter $s$, which depends on a certain design $k\in K_{ad}^H$, can be arbitrarily small when the conductivity of the erzats material $\gamma$ tends to zero. Furthermore, the worst convergence rate takes place at 1-node connected hinges (Fig. \ref{fig:1a}) which form the checkerboard patterns (Example 2.2 in \cite{petzoldt-phd-2001}).

Satisfactory global regularity can be achieved by imposing quasi-monotonicity \cite{dryja-qm-1996} condition on the coefficients.
We say the design $k \in K_{ad}^H$ is \textit{quasi-monotonic at the node} $m \in M^H$ if and only if it has only one local maximum in a small circle around the node $m$ (identifying all maxima lying in the same ground element).
We say the design $k$ is \textit{quasi-monotonic} if and only if it is quasi-monotonic
at each node $m \in M^H$.
Quasi-monotonicity condition is quite a natural restriction in the case of topology optimization problems.
For example, for 0-1 designs \textit{only} 1-node connected hinges 
violate this condition. Similarly, a non-quasi-monotonic node for designs $k \in K_{ad}^H$ is presented in Fig \ref{fig:1b}.
It is important that the quasi-monotonicity condition guarantees that $u(k) \in H^{1+1/4}(\Omega)$ independently of $k$ (Theorem 2.12, Section 2.5 in \cite{petzoldt-phd-2001}).
\begin{figure}[]
    \centering
    \subfloat[]{\label{fig:1a}\includegraphics[width=0.3\columnwidth]{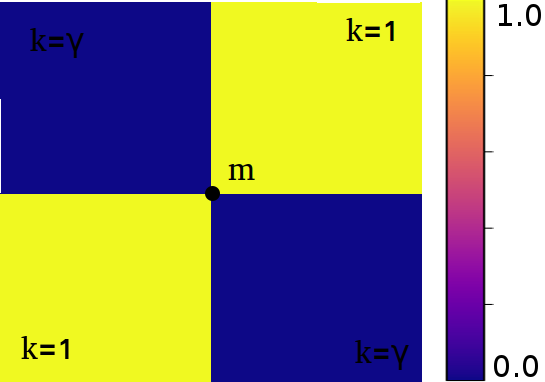}}
    \subfloat[]{\label{fig:1b}\includegraphics[width=0.3\columnwidth]{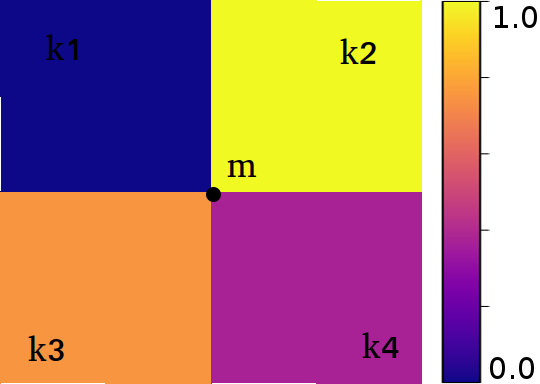}}
    
    \caption{Examples of non quasi-monotonic nodes a) 1-node connected hinge b) analog of 1-node connected hinge for designs admitting intermediate values, assuming $k_1 < k_2,k_3$ and $k_4 < k_2,k_3$.}
    \label{fig:subfigname}
\end{figure}

\end{section}

\begin{section}{A posteriori error estimator}\label{section:apost}
\textit{A priori} error estimate \eqref{convergence}
 describes the asymptotic error behaviour. 
However, we are interested in a technique
that would allow us to \textit{a posteriori} estimate the error for a given finite element solution $u_h$.
In our work, we utilize the estimator considered in \cite{petzoldt-phd-2001}.
It is nothing but a generalization of the estimator for 2D Poisson's equation proposed 
in \cite{verfurth-poissonexplicitapost-1994} to the case of piece-wise constant coefficients.

Let $E_h$ be the set of all edges from $T_h$ and $\omega_E$ be the union of elements that have an edge $E \in E_h$ in common. 
We denote by $k_T$ the value of the coefficient for the element $T \in T_h$
and $k_E = \sum_{T\subset\omega_E}k_T$.
For any \textit{interior} edge $E \in E_h$ and $T,T'\subset \omega_E$, we denote by $n_T$ and $n_{T'}$ the outward normals of $E \subset \partial T$ and $E \subset \partial T'$ respectively.
Given the finite element solution $u_h(k)$, the discrete heat flux $k\nabla u_h$ is a discontinuous across edges $E_h$ function. We define the jump of $k\nabla u_h$  across an interior edge $E$ as follows:
\begin{equation*}
[k\nabla u_h]_E = k_T\dfrac{\partial u_h}{\partial n_T} + k_{T'}\dfrac{\partial u_h}{\partial n_{T'}},
\end{equation*}
assuming the normal derivative operator is applied within corresponding elements.
In the case of homogeneous Neumann boundary condition, the jump of the discrete heat flux across an exterior Neumann edge $E \subset \Gamma_n$ is defined as follows:
\begin{equation*}
 [k\nabla u_h]_E = - k_T\dfrac{\partial u_h}{\partial n_E},
\end{equation*}
where $n_E$ is the outward normal of $E$.

Finally, we present the a posteriori error estimator, denoted as $E_{apost}(k;u_h(k))$, which consists of the sum of local error estimators $\eta_{T}$:
\begin{equation}\label{definition:apost}
\begin{aligned}
 E_{apost}(k;u_h(k)) &= \sum_{T \in T_h}\eta_{T}^2, \\
 \eta_{T}^2 = \dfrac{h^2}{k_T}\|f + \nabla \cdot k \nabla u_h\|^2_{L^2(T)} &+ \sum_{E \subset \partial T/\Gamma_u}\dfrac{h}{k_E}\|[k\nabla u_h]_E\|_{L^2(E)}^2.
\end{aligned}
\end{equation}
It should be noted, that the operator $\nabla \cdot k \nabla (\cdot)$ is also applied within elements, and the term $\nabla \cdot k \nabla u_h$ vanishes for the piece-wise bi-linear solution $u_h(k) \in \mathcal{H}_h^1$, whereas  $E_{apost}(k;u(k)) = 0$ for the true solution $u(k) \in \mathcal{H}$ since the strong residual $f + \nabla \cdot k \nabla u$ vanishes almost everywhere in $\Omega$, and the heat flux $k\nabla u$ is a continuous function.
We also formulate a sufficient condition leading to the robustness of such estimator (Theorem 3.5 in \cite{petzoldt-phd-2001}):
\begin{theorem}\label{theorem:apost}
 If the coefficient $k \in K_{ad}^H$ is quasi-monotonic, 
  then the estimator $E_{apost}(k;u_h(k))$ is robust,
  that is:
 \begin{equation}\label{apost_theorem}
 \begin{aligned}
  \|u(k) - u_h(k)\|_a^2 \leq CE_{apost}(k;u_h(k)),
 \end{aligned}
 \end{equation}
where the constant $C$ only depends on shape regularity of $T_h$.
\end{theorem}
Thus, in the case of quasi-monotonic coefficients, such
 a posteriori error estimator can be regarded as the \textit{robustness indicator}, in the sense that if $E_{apost}(k;u_h(k)) \to 0$ when $h \to 0$ then
$\Phi_h(k) \to \Phi(k)$. 
\end{section}

\begin{section}{Modification of the discrete cost functional using a posteriori error estimator}\label{section:modification}

Most of practically used optimization methods for solving topology optimization problems are based on the sensitivity analysis: following the SIMP approach, they perform continuous gradient-driven minimization of the discrete cost functional and require only the computation of its derivatives with respect to the design variables (so-called sensitivities). Different approaches can be used: the Optimality Criteria methods \cite{bendose-topoptbook-2013}, Sequential Linear Programming (SLP) methods, the Method of Moving Asymptotes (MMA) \cite{svanberg-mma-1987} and other methods for large-scale non-linear constrained optimization.

If we directly minimize the discrete cost functional and do not specifically care about the approximation error, then it is quite natural that optimization algorithms
at each particular gradient step give preference to the ''false minima'' points (designs).

Consider the modification of the discrete cost functional motivated by Theorems \ref{maintheorem}, \ref{theorem:apost}:
 \begin{equation}\label{modfunc}
\Phi_h^{C}(k) = \Phi_h(k) + CE_{apost}(k;u_h(k)),
 \end{equation}
 then, the following apparent corollary holds:
 \begin{corollary}
 Let the design $k \in K_{ad}^H$ be quasi-monotonic, and the constant $C$ is chosen such that \eqref{apost_theorem} holds. Then the modified cost functional $\Phi_h^C(k)$ is the upper bound for the true cost functional $\Phi(k)$.
 \end{corollary}
In order to eliminate the 'false minima' problem, we propose to minimize this modified functional $\Phi^{C}_h(k)$: we do want to minimize $\Phi_h(k)$, but we do not want the FEM error to be too large.
The corresponding modified \eqref{problem:topopt-discrete-fem} problem looks as follows:
\begin{equation}\label{problem:topopt-discrete-fem-modified}\tag{modified $\mathbf{TO}^H_h$}
\left \{
\begin{aligned}
 &\textrm{{\large Minimize}} \;\; \Phi_h^C(k),\\
 (k,&u_h) \in K^H_{ad}\times\mathcal{H}_h\\
  &\textrm{{\large Subject to:}} \\
  &\Phi_h^C(k) = \ell(u_h(k)) + CE_{apost}(k;u_h(k)),\\
   &a_{k}(u_h,v_h)=\ell(v_h),\;\forall v_h \in \mathcal{H}_h.
\end{aligned}
\right.
\end{equation}
Such functional modification can be considered as
a more accurate evaluation of the true cost functional $\Phi(k)$
since it does not affect the solution of the \eqref{problem:topopt-discrete} problem.
It is also worth noting that the constant $C$ from \eqref{apost_theorem} can be found analytically for a certain grid $T_h$.
However, we do not require the designs to be quasi-monotonic during the optimization procedure. We consider $C$ as  the \textit{correction parameter} and investigate the dependence on it in Section \ref{subsec:param_dep}.

We highlight our main observations on solving of the \eqref{problem:topopt-discrete-fem-modified} problem for the considered heat conduction model problem and discretization. Although we do not have theoretical justification of these facts, they are confirmed by our numerical experiments presented in Section \ref{subsec:param_dep}. 
\begin{itemize}
	\item Checkerboard patterns completely disappear already at the small values of the correction parameter $C \approx 0.01$ 
	\item With a certain choices of the parameter $C$, we have managed to get the designs which are very close to the designs obtained using more accurate (and more expensive) finite element approximations
	\item The most important observation is that solving of the \eqref{problem:topopt-discrete-fem-modified} problem with $C$ being large enough gives designs which are very close to the quasi-monotonic ones. This fact is quite surprising and requires further study since the smallness of the estimate does not imply the quasi-monotonicity condition
\end{itemize} It should also be noted that the computation of the sensitivity of $E_{apost}(k;u_h(k))$ requires the solution of the adjoint BVP:
\begin{equation*}
  \begin{aligned}
  \frac{dE_{apost}(k;u_h)}{dk} &= \frac{\partial E_{apost}(k;u_h)}{\partial k} - \lambda^T \frac{\partial A}{\partial k} u_h,\\
  A \lambda &= \frac{\partial E_{apost}(k;u_h)}{\partial u_h},
  \end{aligned}
\end{equation*}
where $A$ denotes the stiffness matrix for the primal BVP.
The complexity in this case nearly doubles, if iterative solvers are used (the same preconditioner can be reused for the adjoint problem).
\end{section}

\begin{section}{Numerical experiments}\label{section:numexp}
\begin{subsection}{Setting up the problem}
We explore the model problem that was previously considered in \cite{gersborg-fvmheat-2006}.
The reference unit square design domain $\Omega$ with the corresponding boundary conditions is presented in Fig. \ref{fig:bc}.
The following set of parameters is used:
\begin{itemize}
 \item The heat source is design-independent and uniform over domain, $f \equiv 10^{-2} \;\text{in}\; \Omega$.
  \item $\gamma = 10^{-3}$ represents conductivity of the erzats material.
 \item The volume constraint is fixed at $\mathbf{V} = 0.4$.
 \item We use the penalization parameter $p=4$, 
 some results for $p = 3$ are also given for the comparison purposes.
 \item We always start the optimization procedure from the uniform distribution of the material.
\end{itemize}
\begin{figure}[]
    \centering
    \includegraphics[width=0.3\columnwidth]{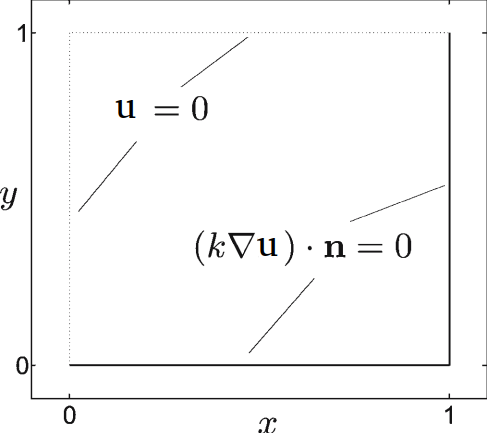}
    \caption{Design domain and boundary conditions.}
    \label{fig:bc}
\end{figure}
We implement the task using several open source software packages:
we use Firedrake package \cite{rathgeber-firedrake-2016} for the finite element analysis and
IPOPT \cite{wachte-ipopt-2006} solver (that implements primal-dual interior-point method) for the optimization.

\end{subsection}

\begin{subsection}{Dependence on the correction parameter}\label{subsec:param_dep}

In this Subsection, we study the dependence of the optimization procedure results 
on the error estimator multiplier $C$ \eqref{modfunc}.
We consider the \eqref{problem:topopt-discrete} problem with the uniform model grids $M^H$, $H=\frac{1}{N}$, $N\in\{64,128\}$.
We solve its \eqref{problem:topopt-discrete-fem-modified} discretization using piece-wise bi-linear approximations on the computational grid that coincides with the model grid ($\mathcal{H}_h = \mathcal{H}_h^1$, $h=H$).
The dependence of the a posteriori error and the cost functional on the correction parameter $C$ is shown in Fig. \ref{fig:apost_fval}, where each resulting design is also evaluated using fine computational grid $T_{h}$, $h=\frac{1}{512}$ for the verification purposes.
We examine the quasi-monotonicity condition using the characteristic function $QM(k)$ presented in the Appendix \ref{sec:appendix}: if $QM(k) = 0$ for the design $k \in K_{ad}^H$ then $k$ does not contain 1-node connected hinges and the quasi-monotonicity condition is satisfied.
The number of iterations until the optimization procedure converges and the quasi-monotonicity values are shown in Fig. \ref{fig:iter_mono}.
Varying the parameter $C$, a lot of qualitatively different designs were obtained,
some of which are presented in Figs. \ref{fig:64alpha_dep} and \ref{fig:128alpha_dep}.

\begin{figure}[]
    \centering
    \subfloat[$N=64$]{\label{fig:64_apost_fval}\includegraphics[width=0.5\columnwidth]{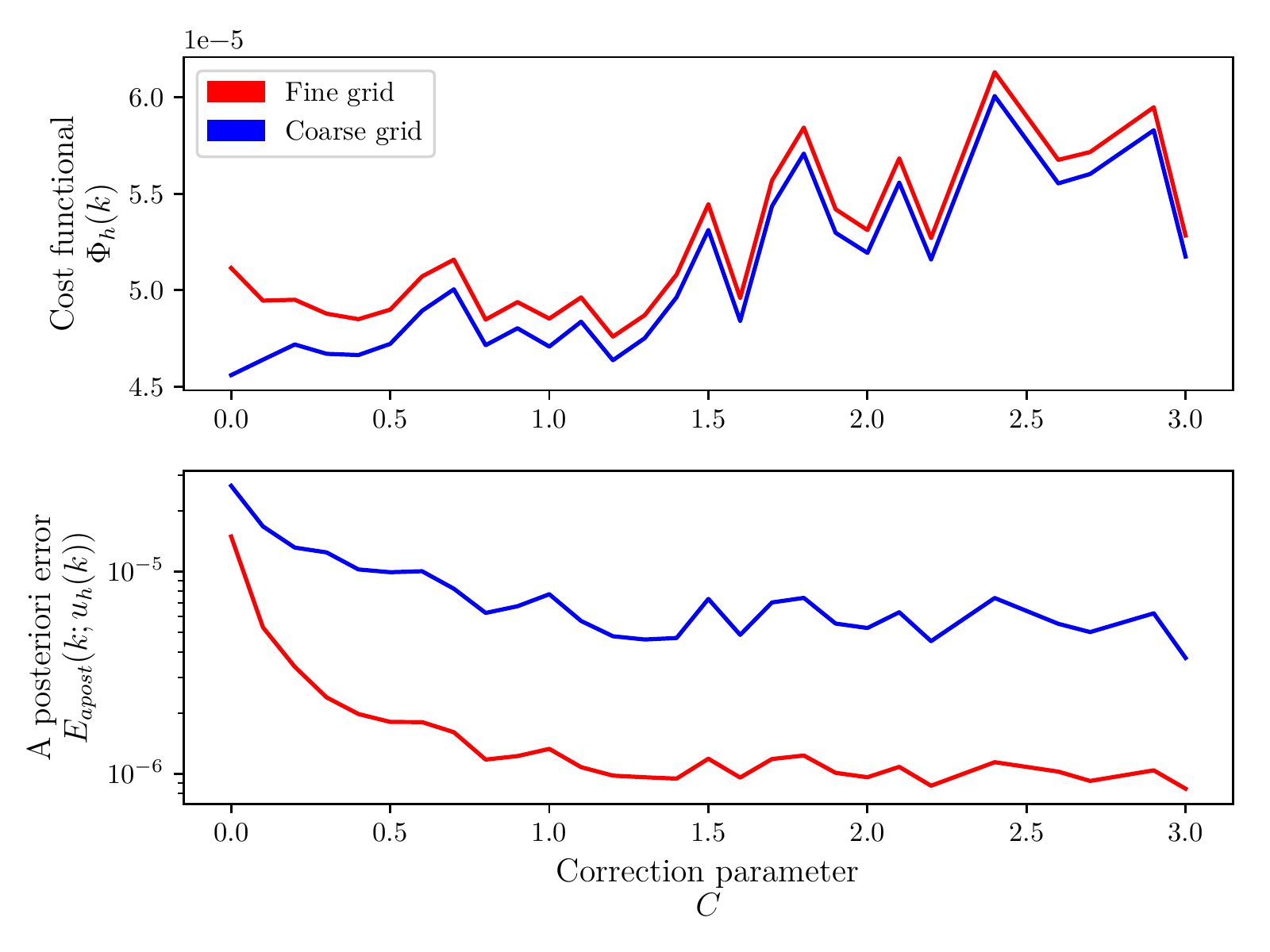}} 
    \subfloat[$N=128$]{\label{fig:128_apost_fval}\includegraphics[width=0.5\columnwidth]{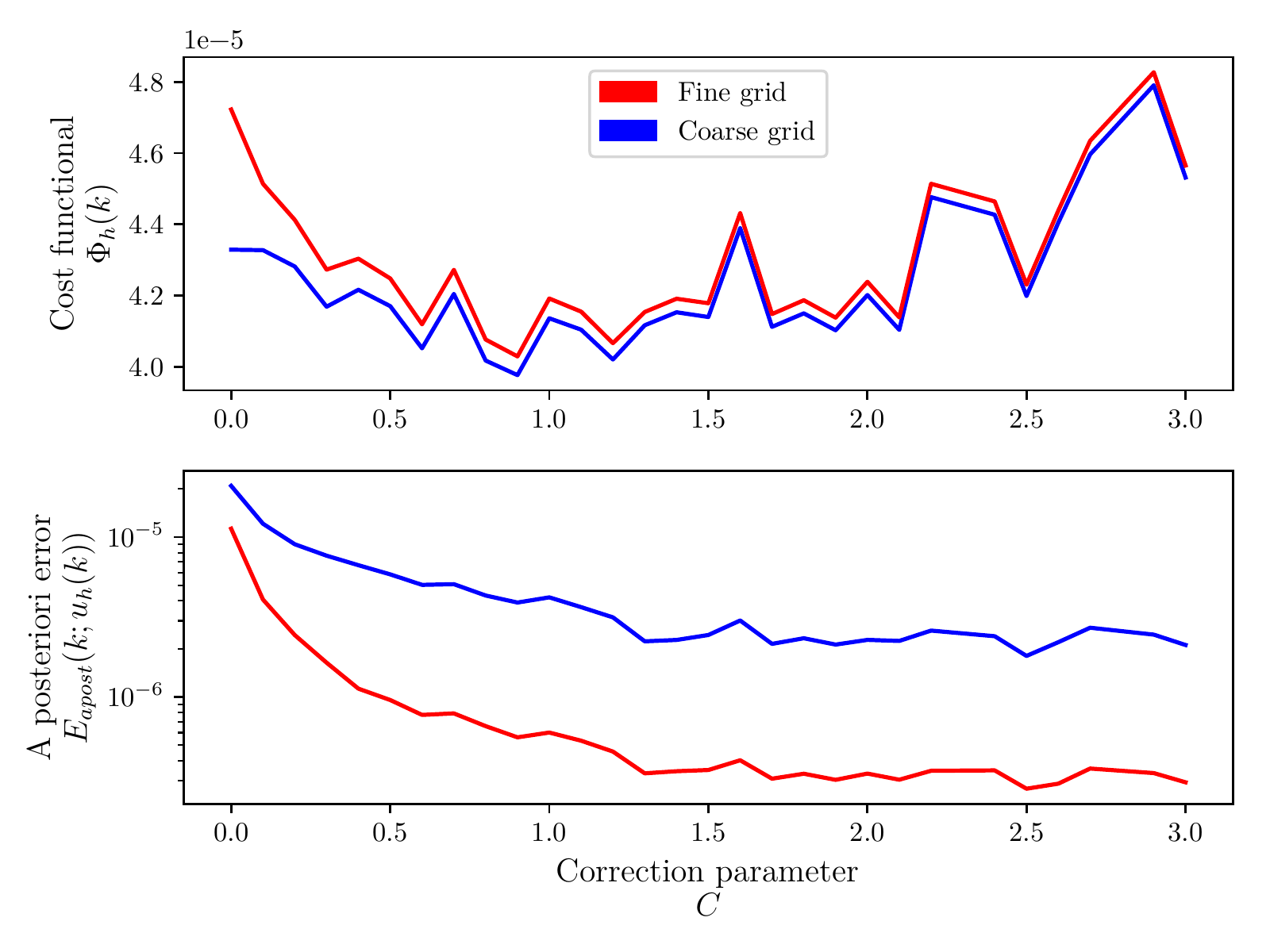}}
    \caption{The dependence of the a posteriori error and the discrete cost functional value on the correction parameter $C$.}
    \label{fig:apost_fval}
\end{figure}

\begin{figure}[]
    \centering
    
    \subfloat[$N=64$]{\label{fig:64_iter_mono}\includegraphics[width=0.5\columnwidth]{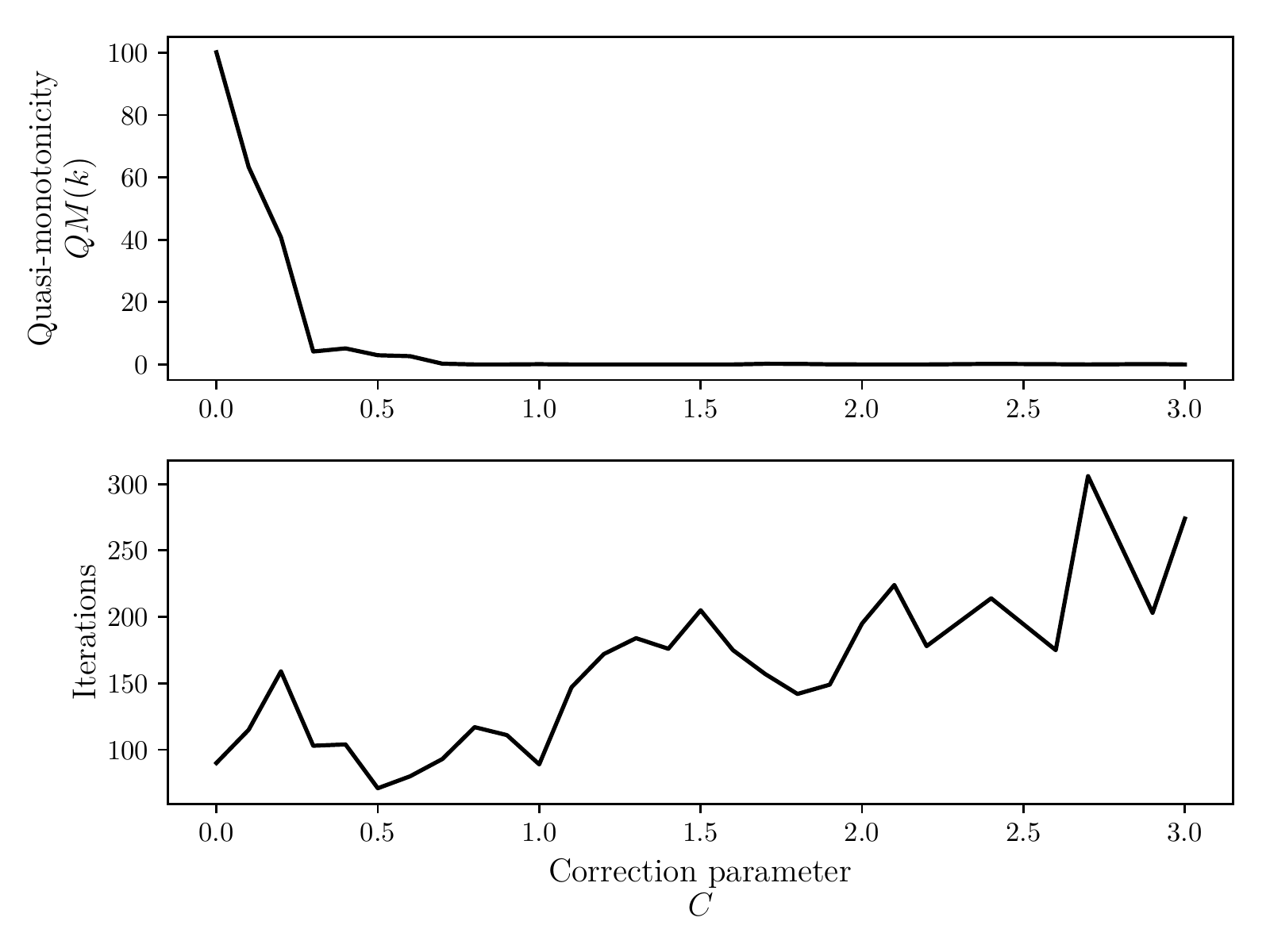}}
    \subfloat[$N=128$]{\label{fig:128_iter_mono}\includegraphics[width=0.5\columnwidth]{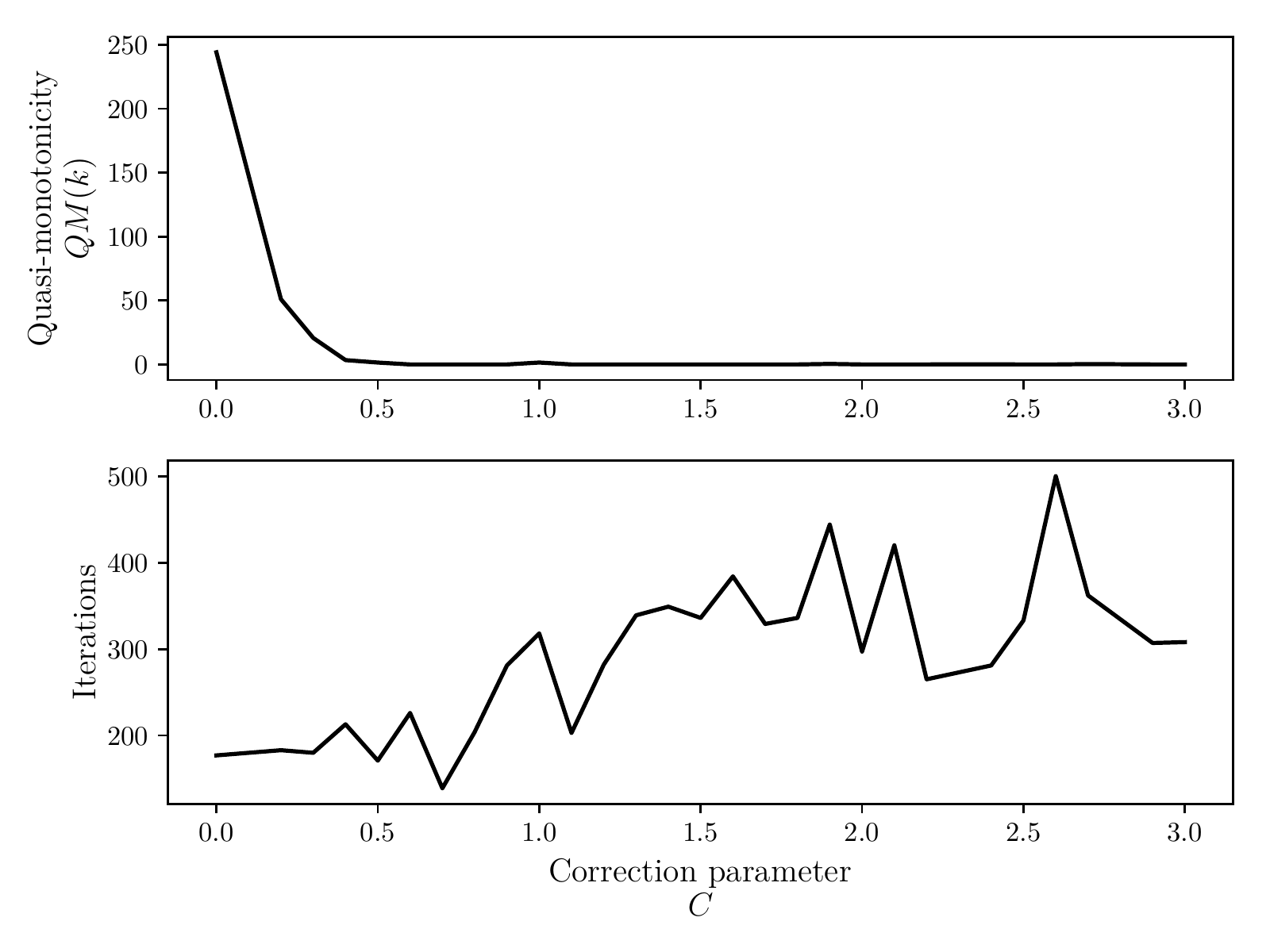}}
    \caption{ The dependence of the number of iterations and the quasi-monotonicity value on the correction parameter $C$ .}
    \label{fig:iter_mono}
\end{figure}

\begin{figure}[]
    \centering
        \subfloat[$C=0.013$]{\label{64013}\includegraphics[width=0.33\columnwidth]{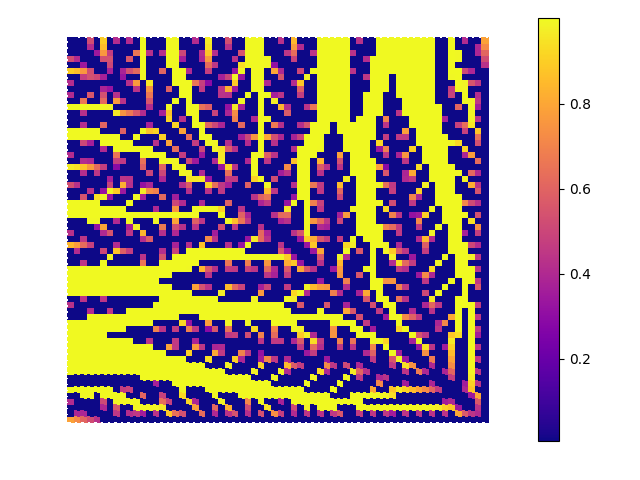}}
        \subfloat[$C=0.2$]{\label{6402}\includegraphics[width=0.33\columnwidth]{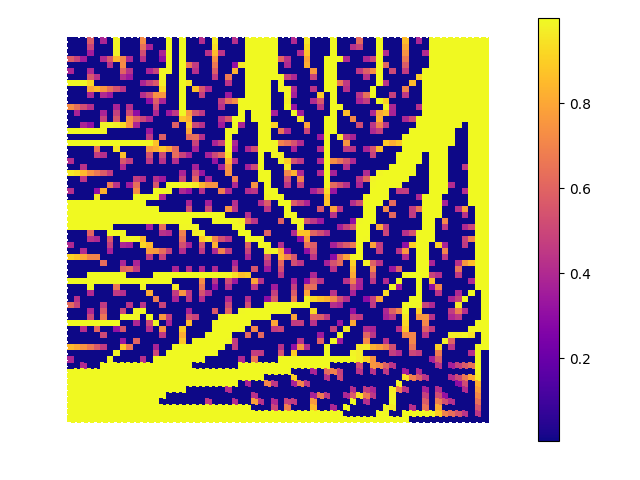}}
        \subfloat[$C=0.4$]{\label{6404}\includegraphics[width=0.33\columnwidth]{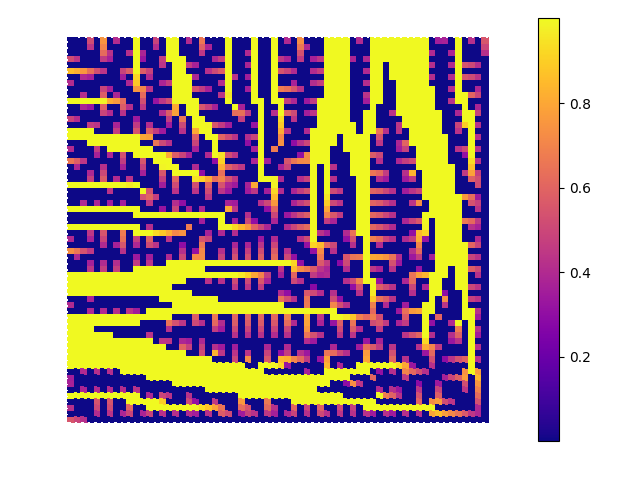}}
        
        \subfloat[$C=0.8$]{\label{6408}\includegraphics[width=0.33\columnwidth]{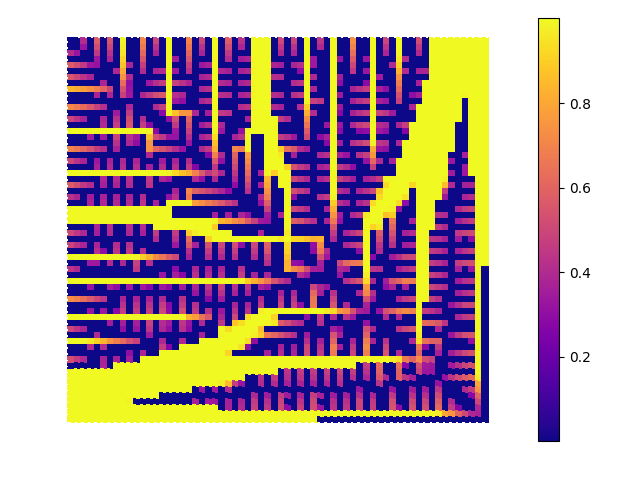}}
        \subfloat[$C=3.0$]{\label{6430}\includegraphics[width=0.33\columnwidth]{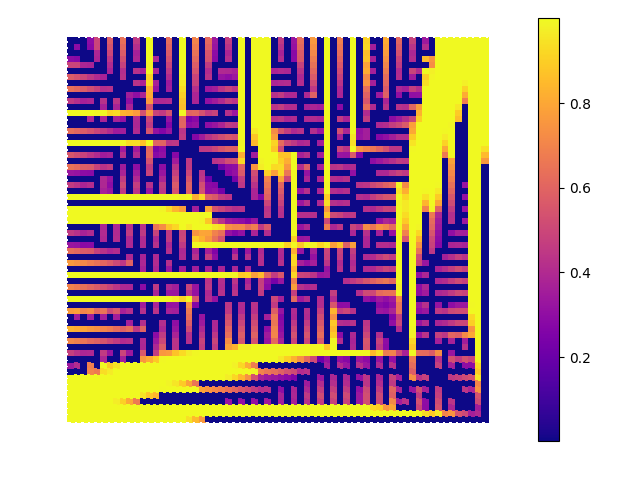}}
    \caption{Designs for different correction parameters $C$, $N=64$}
    \label{fig:64alpha_dep}
\end{figure}

\begin{figure}[]
    \centering
        \subfloat[$C=0.0$]{\label{12800}\includegraphics[width=0.3\columnwidth]{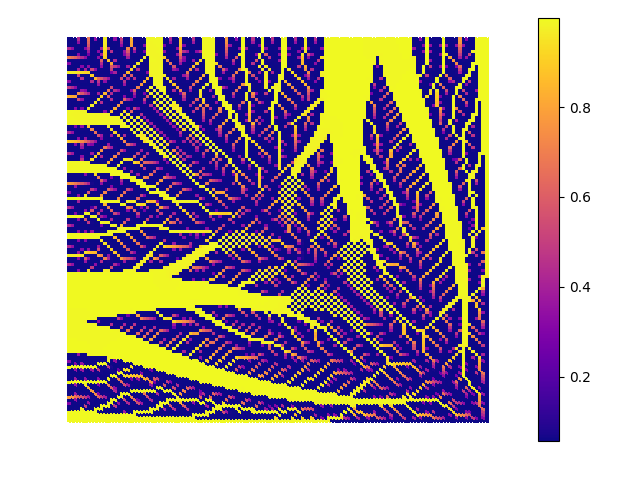}}
        \subfloat[$C=0.1$]{\label{12801}\includegraphics[width=0.3\columnwidth]{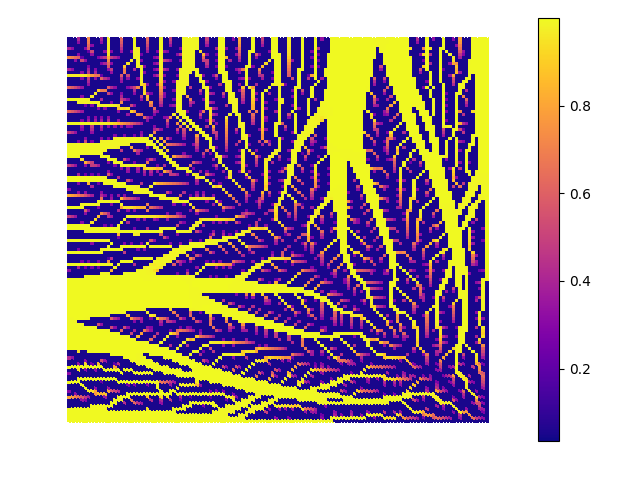}}
        \subfloat[$C=0.3$]{\label{12803}\includegraphics[width=0.3\columnwidth]{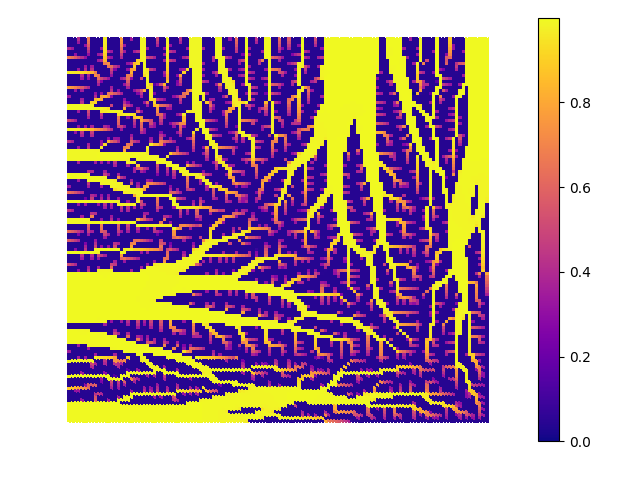}}
        
        \subfloat[$C=0.4$]{\label{12809}\includegraphics[width=0.3\columnwidth]{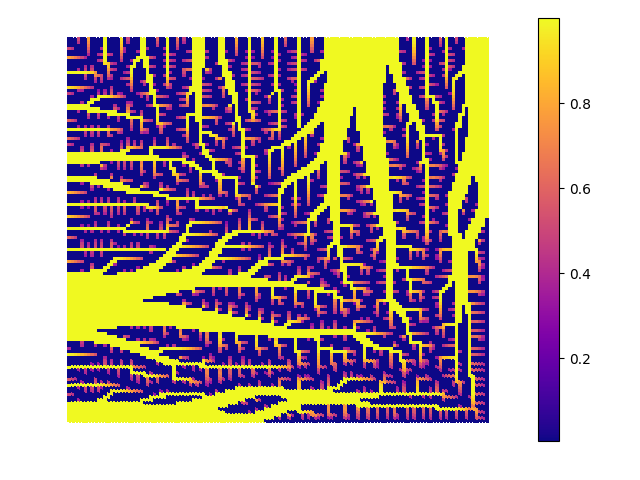}}
        \subfloat[$C=0.9$]{\label{12830}\includegraphics[width=0.3\columnwidth]{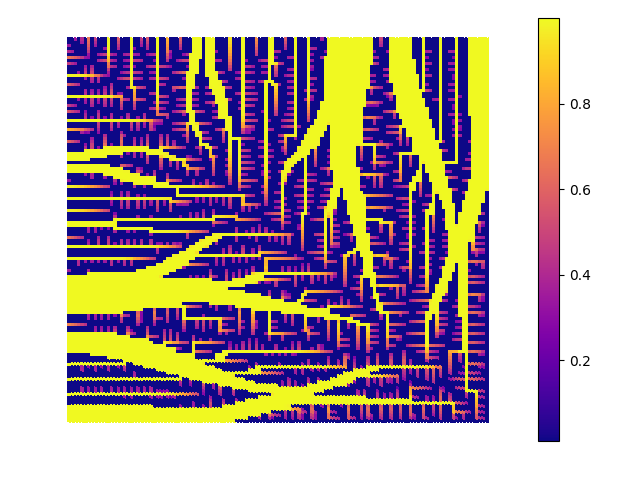}}
    \caption{Designs for different correction parameters $C$, $N=128$.}
    \label{fig:128alpha_dep}
\end{figure}

Below we highlight the main observations that are valid for the considered problem and discretization.
Although we do not have theoretical justification of these facts,
they are confirmed by our numerical experiments:
\begin{itemize}
\item The error decreases when increasing the correction parameter until $ C \approx 1.5$
\item The decrease in the error naturally leads to the eventual elimination of the checkerboard patterns when $C \gtrapprox 0.09$
\item Following Fig. \ref{fig:iter_mono}, we can observe that the quasi-monotonicity condition is automatically provided when $C \gtrapprox 0.6$ that implies the values obtained using the fine computational grid are robust. As mentioned above, this fact is quite surprising since whereas it is true that the quasi-monotonicity condition implies the robustness of the estimate, it should be emphasized that the smallness of the a posteriori estimate does not imply the quasi-monotonicity condition
\item Another important result is that, with certain choices of the parameter C, we have managed to get the designs which are very close (both visually and by the value of the cost functional) to the designs obtained using more accurate (and more expensive) approximations.
For example, the design in Fig. \ref{64013} is very close to the design obtained using once-refined
computational grid (Fig. \ref{1ref}), when the design in Fig. \ref{6402} is very close to the design obtained using bi-quadratic finite element approximation (Fig. \ref{2lag})
\end{itemize}

The correction implies the increase in the regularity parameter $s$ \eqref{convergence} since the designs become quasi-monotonic.
The decrease in the error can be also explained by the straightening of the “streaks” of the designs, since the number of singular nodes (corners) in the design $k$ directly determines the term $|u(k)|^2_{H^{1+s}}$.


The choice of the optimal correction parameter is not entirely clear since even small changes in the parameter can lead to the falling into different local minima. However, we are encouraged by the fact that the values of the parameter at which the checkerboards disappear and the quasi-monotonicity condition is provided are very close for the different model sizes. Based on our experiments, we can conclude that the choice $C \approx 1.0$ is appropriate for the considered model problem.

In Fig. \ref{fig:convergence10} we compare the convergence of the optimization procedure for the fixed parameters $C = 0.0$  and $C = 1.0$. It is clearly seen how the correction helps to suppress the error in the later case.  
We also provide the designs for some intermediate iterations of the optimization process (Fig. \ref{fig:convergence_steps}).

\begin{figure}[]
    \centering
    \includegraphics[width=0.6\columnwidth]{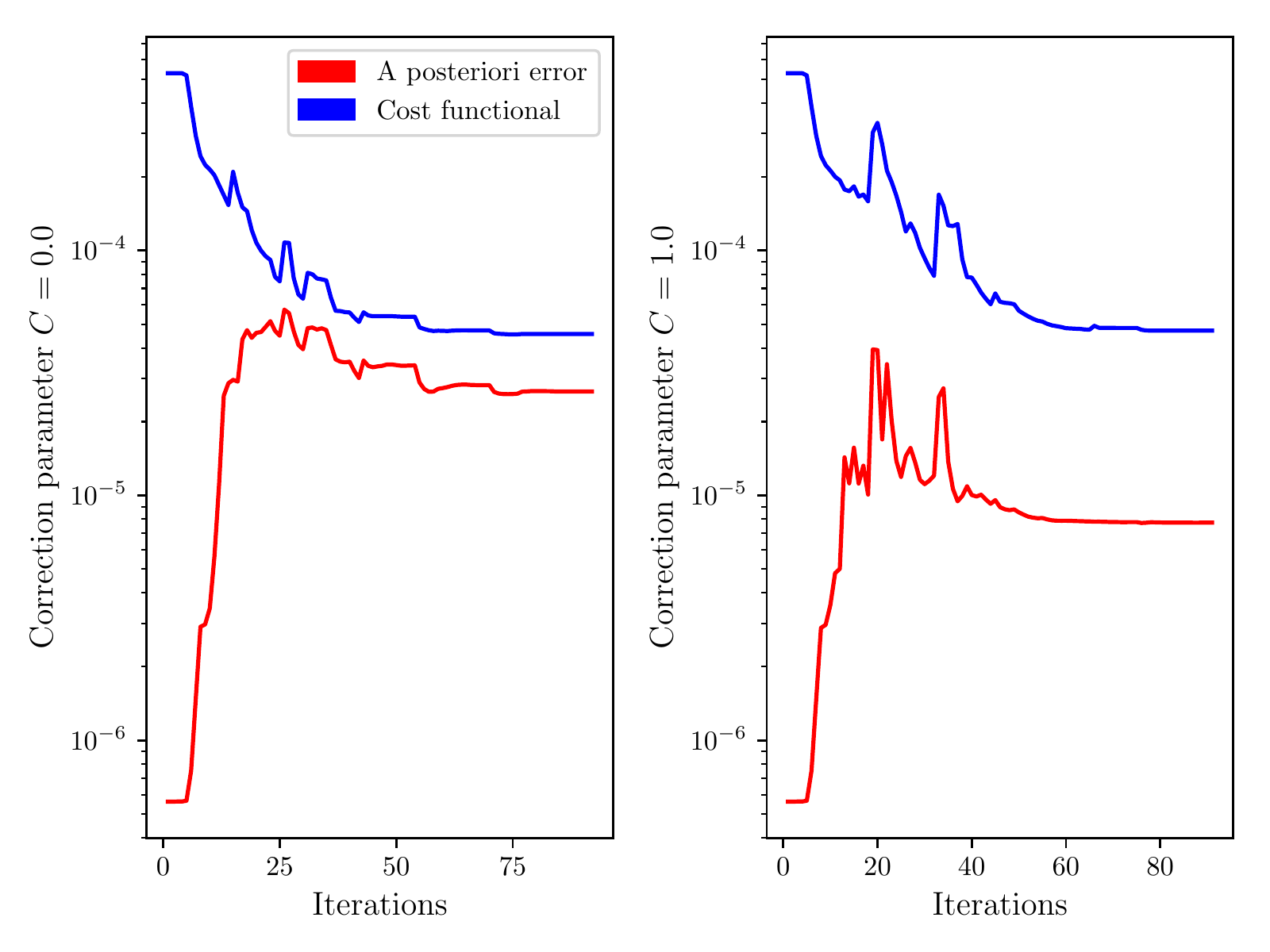}
    \caption{The convergence of the cost functional and the a posteriori error. $N=64$.}
    \label{fig:convergence10}
\end{figure}

\end{subsection}

\begin{subsection}{Demonstration of the ``false minima'' problem}\label{subsec:appr_problem}
All presented in Fig. \ref{fig:designs} designs are the solutions of the \eqref{problem:topopt-discrete} problem with the fixed model grid $M^H$, $N=64$ obtained using different approximations:
\begin{figure}[]
    \centering
        \subfloat[]{
        \label{cb}
        \includegraphics[width=0.33\columnwidth]{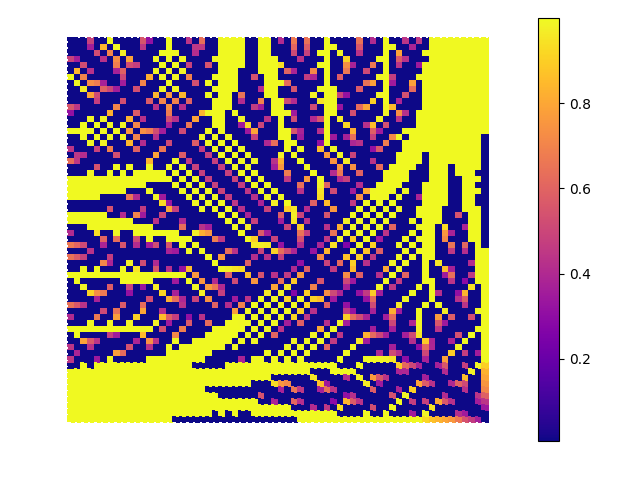}}
        \subfloat[]{
        \label{1ref}
        \includegraphics[width=0.33\columnwidth]{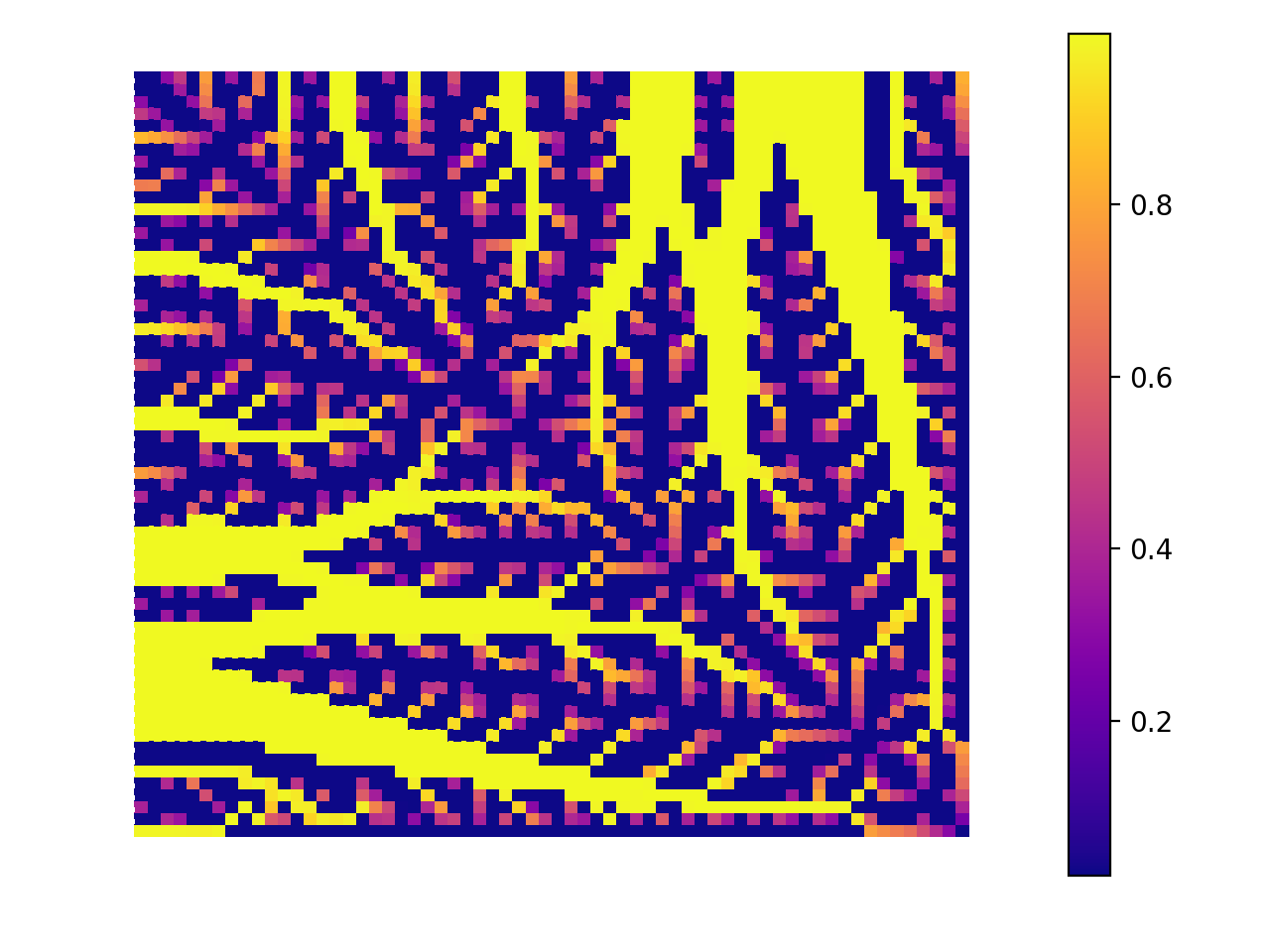}}
        \subfloat[]{
        \label{2lag}
        \includegraphics[width=0.33\columnwidth]{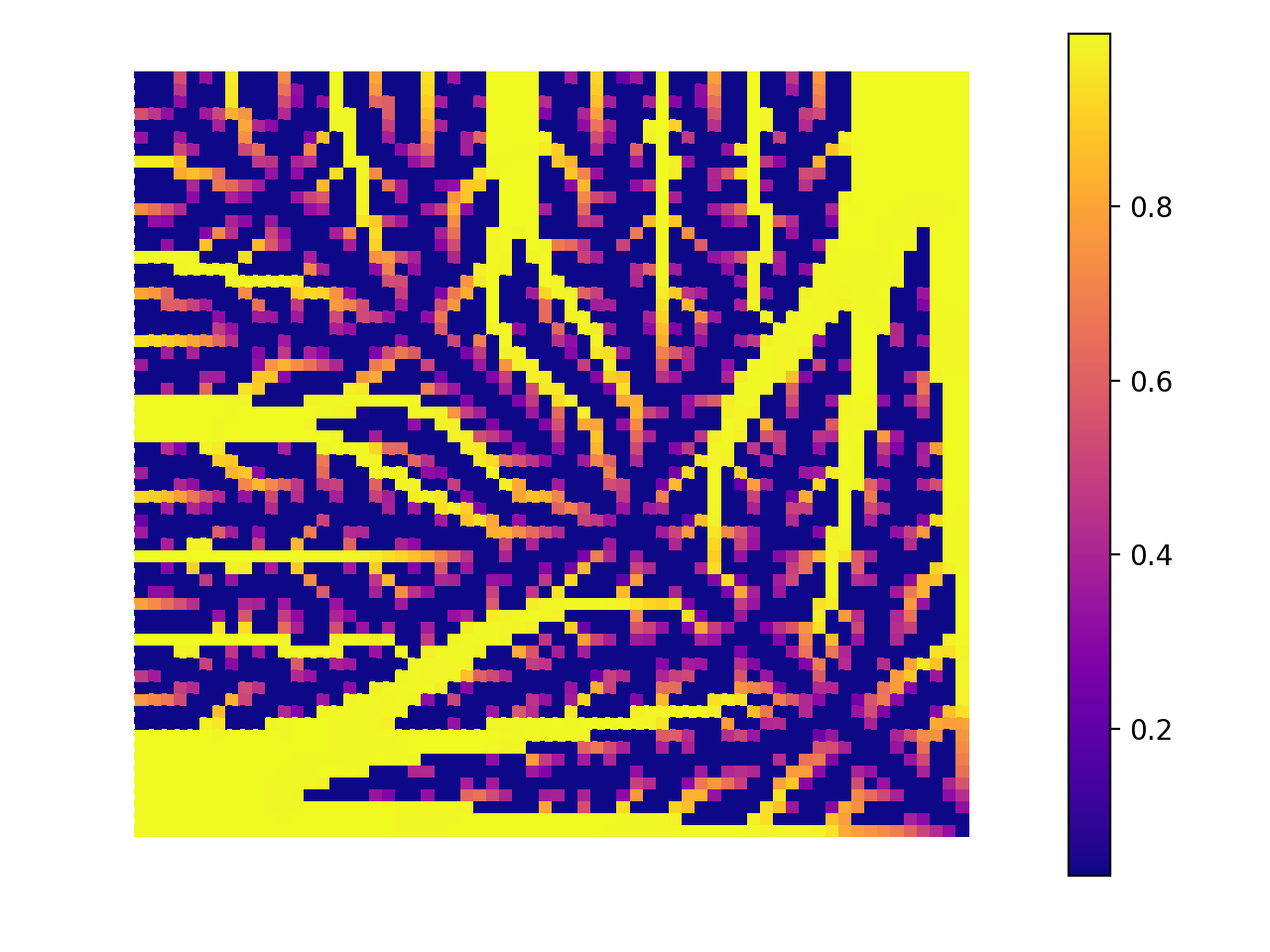}}
        
        \subfloat[]{
        \label{filter1}
        \includegraphics[width=0.33\columnwidth]{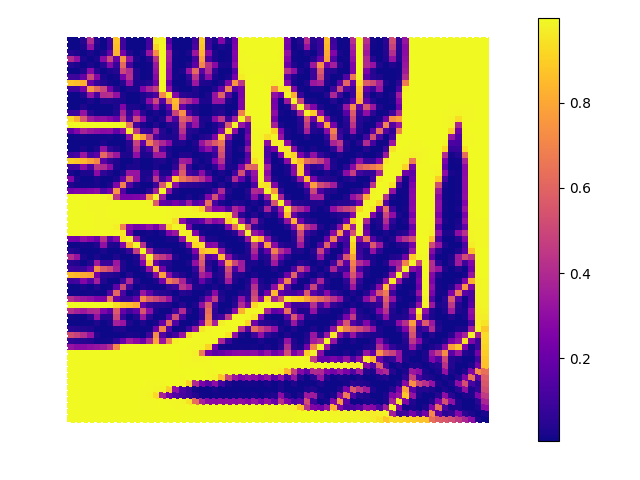}}
        \subfloat[]{
        \label{filter2}
        \includegraphics[width=0.33\columnwidth]{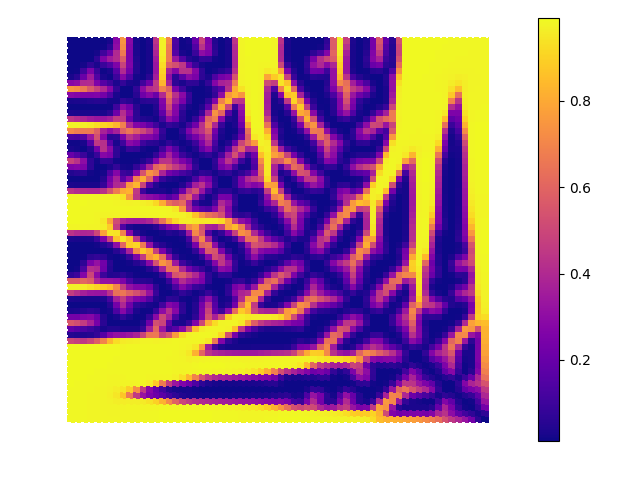}}
        \subfloat[]{
        \label{our}
        \includegraphics[width=0.33\columnwidth]{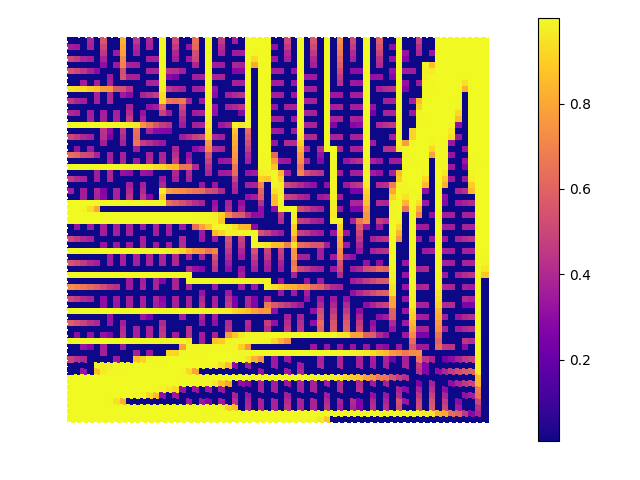}}
    \caption{a) typical checkerboard b) once-refined computational grid c) bi-quadratic finite elements
    d) small-radius sensitivity filter e) large-radius sensitivity filter
    f) proposed functional modification, $C = 1.2$}
    \label{fig:designs}
\end{figure}

\begin{figure}[]
    \centering
    \includegraphics[width=0.8\columnwidth]{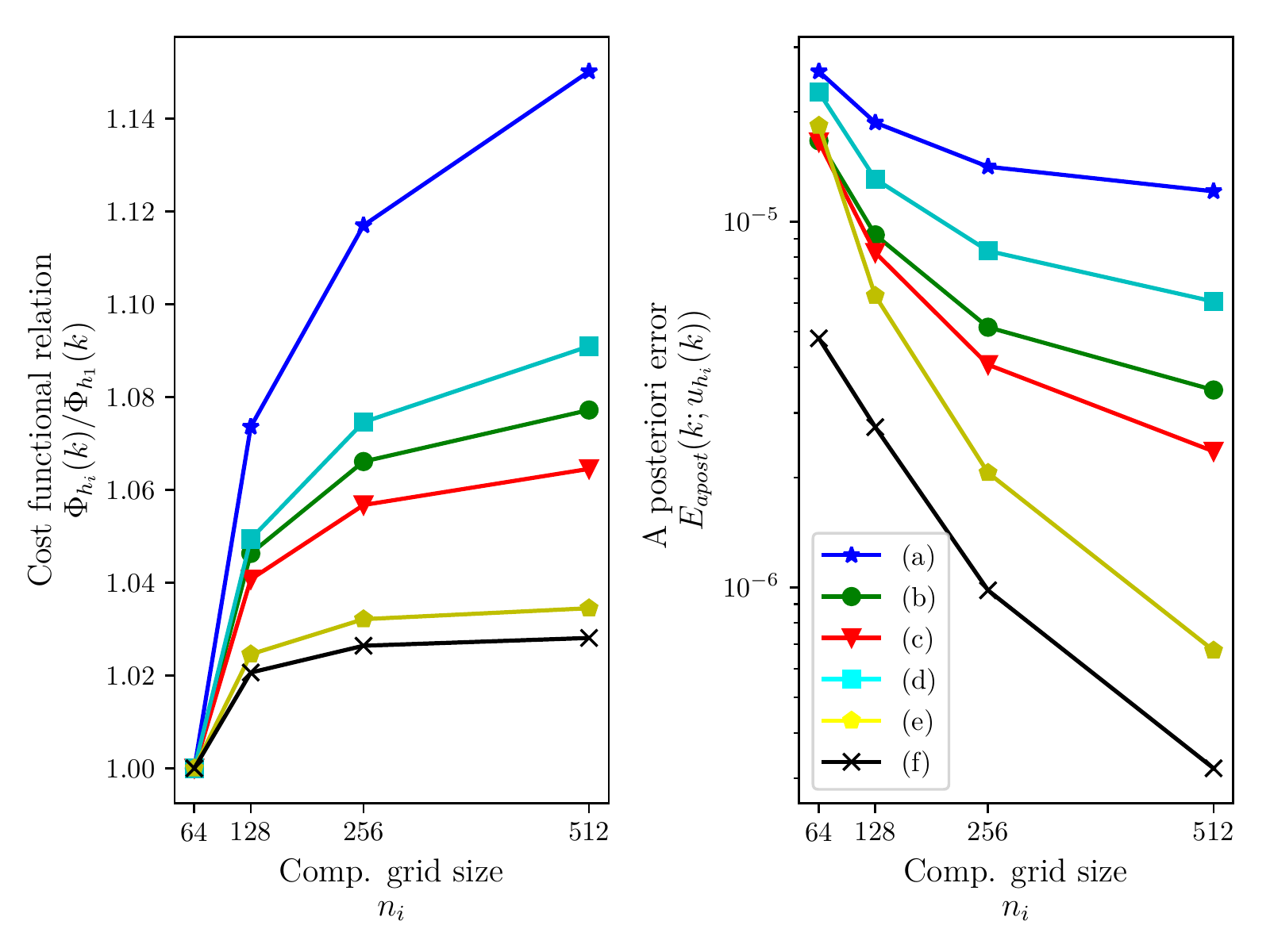}
    \caption{Refinement results: a) typical checkerboard b) once-refined computational grid c) bi-quadratic finite elements d) small-radius sensitivity filter e) large-radius sensitivity filter
    f) proposed functional modification, $C = 1.2$}
    \label{fig:refinement_results}
\end{figure}
\begin{itemize}

 \item Fig. \ref{cb} shows the design where the checkerboard problem is clearly traced.
It is obtained solving the \eqref{problem:topopt-discrete-fem} problem with the standard piece-wise bi-linear approximation on the computational grid that coincides with the model grid ($\mathcal{H}_h = \mathcal{H}_h^1$, $h=H$)

\item Figs. \ref{1ref} and \ref{2lag} show the designs
obtained solving \eqref{problem:topopt-discrete-fem} problem
using once-refined computational grid ($\mathcal{H}_h = \mathcal{H}_h^1$, $2h=H$)
and bi-quadratic finite element approximation ($\mathcal{H}_h = \mathcal{H}_h^2$, $h=H$) respectively

\item We also consider the designs obtained using the classical sensitivity filter \cite{sigmund-99line-2001,bendose-topoptbook-2013}
that is widely used to prevent both the checkerboard problem and the mesh-dependency phenomena.
In Fig. \ref{filter1} the filter radius is chosen so that only checkerboards
are removed, and in Fig. \ref{filter2} the radius is large enough to provide mesh-independency

\item Finally, we present the design obtained solving the \eqref{problem:topopt-discrete-fem-modified} problem with $\mathcal{H}_h = \mathcal{H}_h^1$ and $h=H$.
The optimal parameter $C = 1.2$ is chosen by inspecting the dependence shown in Fig. \ref{fig:64_apost_fval}.
The design resembles a lamellar needle structure that correlates with the work \cite{yan-nonoptimalityoftrees-2018} on the non-optimality of tree-like structures for heat conduction problems.
\end{itemize}

For each resulting design we investigate the behaviour of the discrete cost functional and the a posteriori error
on the family of conformal computational grids
$\{T_{h_i}\}, \; h_i = \frac{1}{n_i}$,
 $n_i \in \{64, 128, 256, 512\}$,
where the coarsest one coincides with the model grid $M^H$.
On each grid we use piece-wise bi-linear approximations $\mathcal{H}_h^1$.
The corresponding refinement results are shown in Fig. \ref{fig:refinement_results}:
the dependence of the relation $\Phi_{h_i}(k)/\Phi_{h_1}(k)$ on the computational grid size $n_i$
from the left, and the dependence of the a posteriori error $E_{apost}(k;u_{h_i}(k))$ from the right.
The error and the cost functional values on the finest and the coarsest meshes as well as the quasi-monotonicity value $QM(k)$ are presented in Table \ref{table:1}.

As it was expected for the checkerboard-like design \ref{cb}, we can observe the significant increase in the cost functional when refining the computational grid.
The improved approximations (\ref{1ref}, \ref{2lag}) help to avoid checkerboard patterns, however they do not completely solve the ``false minima'' problem
since even after 8 times refinement,
the error is still of the same order as the functional, and the values $\Phi_{h_4}(k)$ are expected to increase further.
Moreover, we can not sharply estimate
their upper bounds since they are not quasi-monotonic and the a posteriori error estimates are not robust.
The design in Fig. \ref{filter2} is quasi-monotonic so the error estimate is robust.
However, smoothing effect of the filter imposes a very significant restriction
on the original set of admissible designs that entails
a much larger value of the cost functional.
Finally, the proposed modification also provides
the quasi-monotonicity and the small error together with the small and robust cost functional value.

\begin{table}[]
\begin{center}
\begin{tabular}{ c | c | c| c | c | c | c  }

        & (a)          & (b)     & (c) & (d) & (e) & (f) \\ 
\hline
$\Phi_{h_1}(k)$ &  4.57e-05 & 4.46e-05 &4.72e-05 & 5.04e-05& 6.24e-05 & 4.63e-05\\ 
$\Phi_{h_4}(k)$ & 5.25e-05 & 4.80e-05 &5.03e-05& 5.49e-05 & 6.45e-05 & 4.76e-05\\ 
$E_{apost}(u_{h_1}(k))$ & 2.57e-05 & 1.66e-05 & 1.65e-05& 2.25e-05 & 1.83e-05 & 4.79e-06\\ 
$E_{apost}(u_{h_4}(k))$ & 1.21e-05 & 3.46e-06 & 2.35e-06 & 6.04e-06 & 6.72e-07 & 3.20e-07\\ 
$QM(k)$ & 107.06 & 82.21 & 50.04 & 7.13 & 3e-3 & 2e-6\\ 
\end{tabular}

\caption{Refinement results: a) typical checkerboard b) once-refined computational grid c) bi-quadratic finite elements d) small-radius sensitivity filter
e) large-radius sensitivity filter
    f) proposed functional modification, $C = 1.2$}
\label{table:1}
\end{center}
\end{table}

\end{subsection}

\begin{subsection}{Model grid refinement study}

In this Subsection, we present the model grid refinement study for the fixed correction parameter $C = 1.0$. We consider the \eqref{problem:topopt-discrete} problems using the family of model grids
$\{M_{H_i}\}$, $i \in  \{1,\ldots,8\}$, $N_i = 32 \cdot i$. 
We solve the \eqref{problem:topopt-discrete-fem-modified} problems using piece-wise bilinear approximations $\mathcal{H}_{h_i} = \mathcal{H}_{h_i}^1$, $h_i = H_i$.
The dependence of the cost functional and the a posteriori error on the model grid size is shown in Fig. \ref{fig:model_refinement}.
 Some resulting designs are presented in Fig. \ref{fig:model_refinement_designs}.
 The results are mesh-dependent. 
 All obtained designs are quasi-monotonic.
 As it was expected, when refining the model grid, we can observe the decrease in the cost functional, although it is stabilized at the value $3.75\times 10^{-5}$, and the results do not improve when $N_i > 160$. 
 
\begin{figure}[]
    \centering
    \includegraphics[width=0.8\columnwidth]{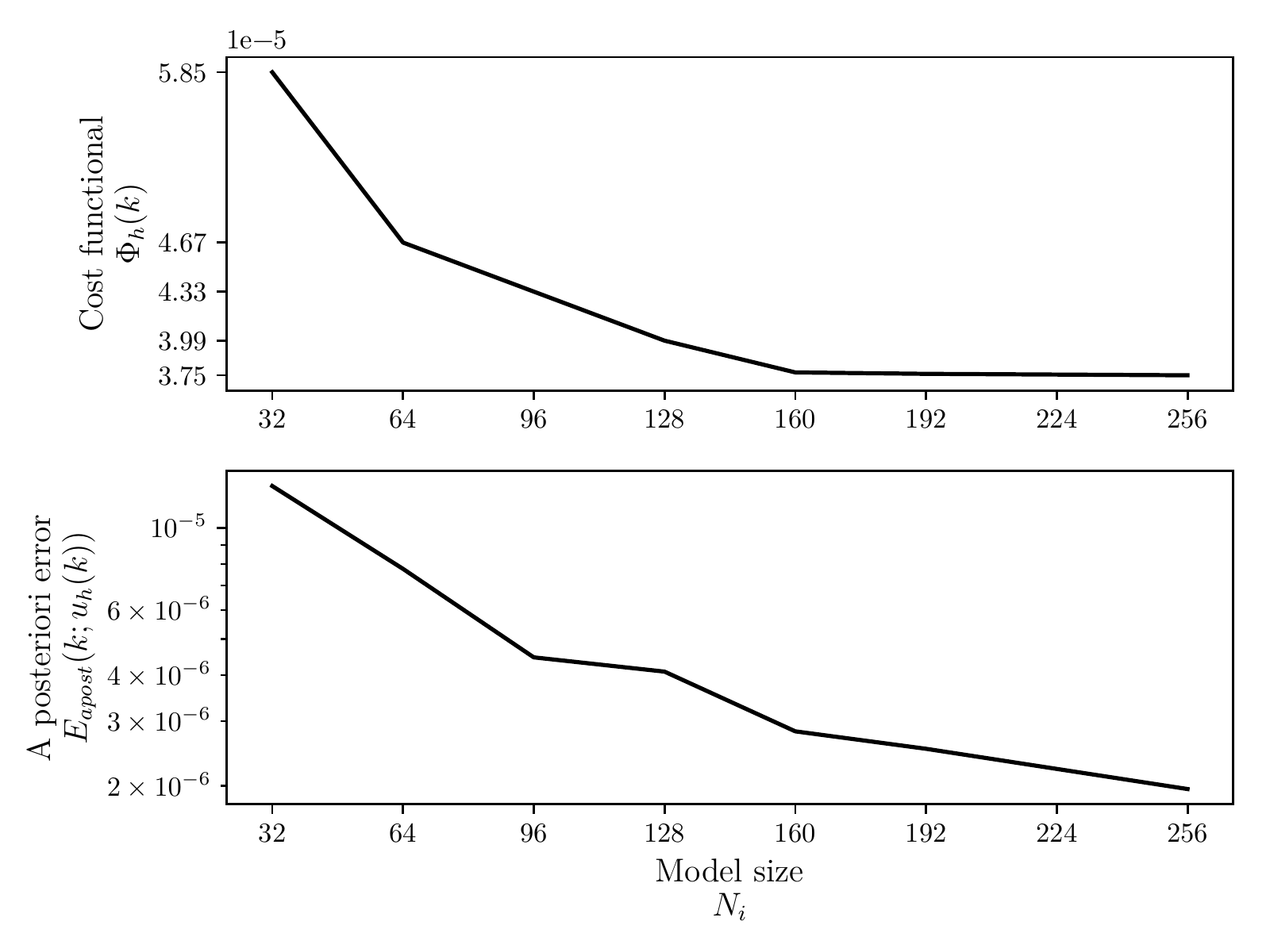}
    \caption{The dependence of cost functional and a posteriori error on model grid size, $C=1.0$}
    \label{fig:model_refinement}
\end{figure}

\end{subsection}

\begin{subsection}{Comparison of results}
In this Subsection, we show some designs for the penalization parameter $p=3$ (Fig. \ref{fig:optimal_p3}). The discrete cost functional value and the a posteriori error estimate computed on the fine computational grid ($n=512$) using piece-wise bi-linear approximations $\mathcal{H}_h=\mathcal{H}_h^1$ can be found in Table \ref{table:2}.
The designs are quasi-monotonic hence the presented values are robust.
It should be noted that, even although the continuation approach \cite{bendose-topoptbook-2013} was not adopted in our study, we were able to get better value of the cost functional than it was reported in \cite{gersborg-fvmheat-2006} where the design with the value $3.82\times10^{-5}$ for $N=128$ was given.

\begin{figure}[]
    \centering
        \subfloat[$C=1.1, N=64$]{
        \label{64p3}
        \includegraphics[width=0.33\columnwidth]{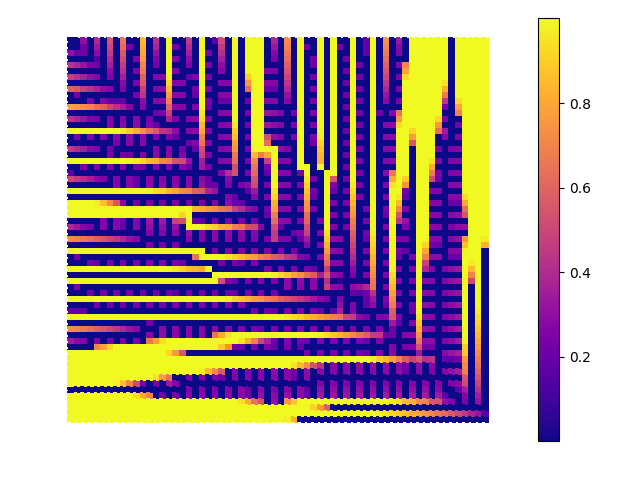}}
        \subfloat[$C=0.8, N=128$]{
        \label{128p3}
        \includegraphics[width=0.33\columnwidth]{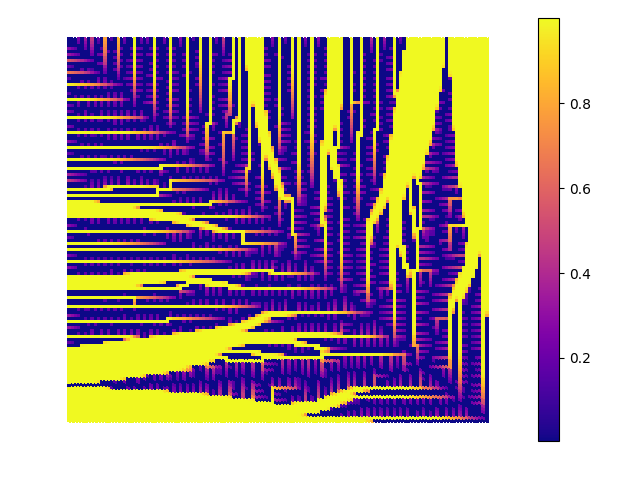}}
    \caption{The best obtained designs for the penalization parameter $p=3$.}
    \label{fig:optimal_p3}
\end{figure}

\begin{table}[h!]
\begin{center}
\begin{tabular}{ c | c| c }
        & Fig. \ref{64p3}     & Fig. \ref{128p3}  \\ 
\hline
$\Phi_{h}(k)$ & 3.94e-05 & 3.62e-05\\ 
$E_{apost}(u_{h}(k))$ & 3.29e-07 & 4.08e-07\\  
\end{tabular}
\end{center}

\caption{Penalization parameter $p=3$.
Robust cost functional and a posteriori error estimate computed on the fine mesh
for the designs presented in Fig. \ref{fig:optimal_p3}.}
\label{table:2}
\end{table}	

\end{subsection}

\end{section}

\begin{section}{Discussion and related works}
%

 Despite the simplicity of the considered heat conduction model problem, the idea can be extended to the other objectives as well as other elliptic problems.  
 The proposed functional modification can be incorporated with any approach whenever the standard piece-wise polynomial finite element approximations are used. 
It is quite interesting whether this idea can give any significant advantages in solving of more complicated three-dimensional problems.
 It should be also clarified, that although the a posteriori error estimator $E_{apost}(k;u_h(k))$ \eqref{definition:apost} contains a non-differentiable jump operator, it is actually smoothed due to the squaring.
 The level set method \cite{vanDijk-levelsetreview-2013,sethian-levelset-2000,wang-levelset-2003,allaire-levelset-2004} should also be mentioned, since most of the level set based approaches also operate with erzats materials and fixed domains and can be considered as a modification of the classical density-based approach \cite{sigmund-lastreview-2013}.
A higher accuracy in the case of discontinuous coefficients
also can be achieved by improving the approximation properties of the standard piecewise polynomial
finite element spaces by enriching it with special functions that 
better approximate a priori known local singularities of the solutions.
For example, the extended finite element method (X-FEM) along with the level set approach was considered
in \cite{wei-xfemlevelset2010}, whereby considerably more accurate results around the interfaces were achieved.
The advantages of the X-FEM were also demonstrated in \cite{guo-mmc-2014} together with
a new Moving Morphable Components (MMC) based framework \cite{guo-mmc-2016,zhang-mmc-2016}.
The implementation of the ESO (Evolutionary Structural Optimization) method \cite{xie-simpleeso-1993,huang-esoreview-2010,munk-evoreview-2015} for steady heat conduction was presented in \cite{li-esoheat-1999},
another related work with a similar model problem is \cite{gao-besoheat-2008}, where the BESO (Bi-directional ESO) method \cite{querin-besoorigin-1998,querin-besovalidation-2000} was used.
Also, the implementations of the level-set method for the heat conduction problems were considered in
\cite{zhuang-heatlevelset-2007,coffin-heatlevelset-2016}.
\end{section}

\appendix
\section{Quasi-monotonicity characteristic function}
\label{sec:appendix}
A simple scheme providing the quasi-monotonicity condition to prevent checkerboards and 1-node connected hinges was considered in \cite{poulsen-qm-2002}, where the characteristic function that detects non-quasi-monotonic designs was presented.
In the case of square design domain $\Omega$, assuming the model grid $M^H$ consists of $N \times N$ ground elements, this characteristic function looks as follows:
\begin{equation*}\label{qm}
 QM(k) = \sum_{j=1}^{N-1}\sum_{i=1}^{N-1} qm(k_{i,j},k_{i+1,j},k_{i,j+1},k_{i+1,j+1}),
\end{equation*}
where $k_{i,j}$ denotes the value of the coefficient in the corresponding ground element $M_{i,j} \in M^H$, and local function $qm$ is a function of four ground elements surrounding a node in the interior of the design:
\begin{equation*}
\begin{aligned}
 qm(a,b,c,d) &= m(a,b,c)\cdot m(a,c,d) \cdot m(b,a,c) \cdot m(b,d,c),\\
 m(a,b,c) &= |b - a| + |c - b| - |c - a|.
\end{aligned}
\end{equation*}

\begin{figure}[]
    \centering
        \subfloat[$N=32$]{\label{mesh:32}\includegraphics[width=0.3\columnwidth]{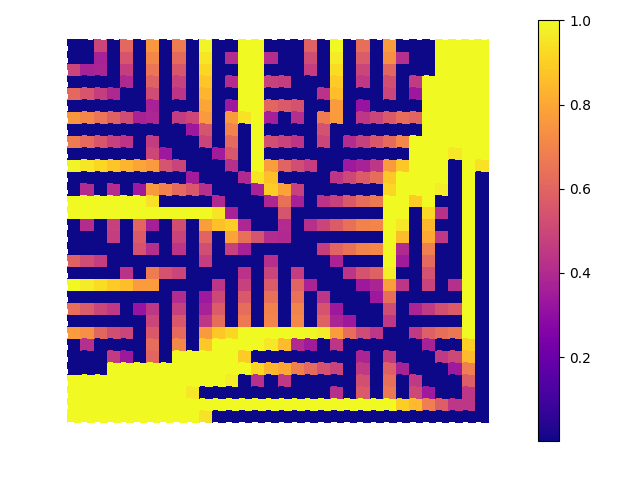}}
        \subfloat[$N=64$]{\label{mesh:64}\includegraphics[width=0.3\columnwidth]{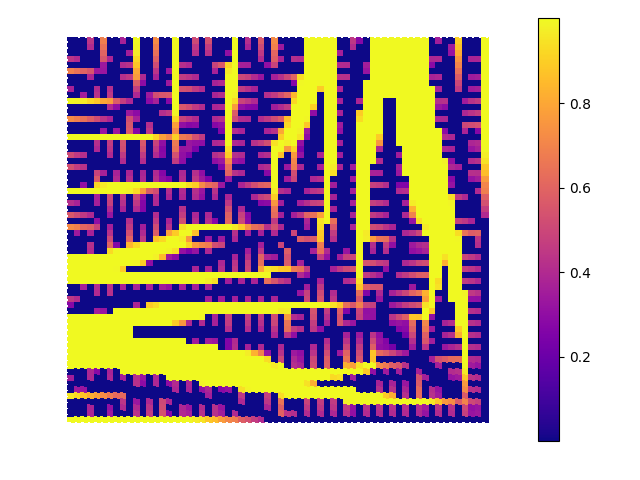}}
        
        \subfloat[$N=96$]{\label{mesh:96}\includegraphics[width=0.3\columnwidth]{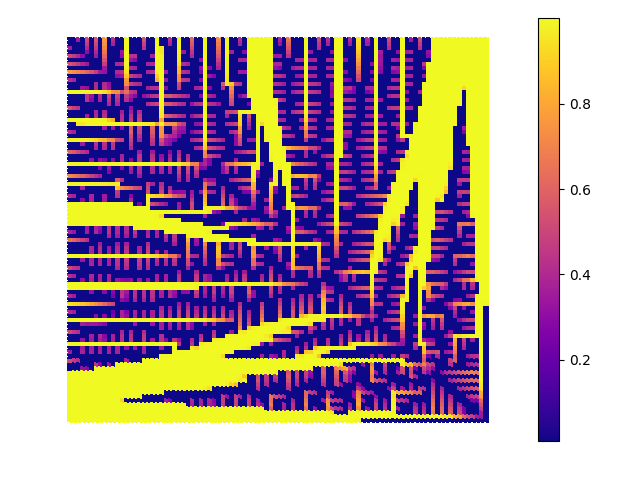}}
        \subfloat[$N=128$]{\label{mesh:128}\includegraphics[width=0.3\columnwidth]{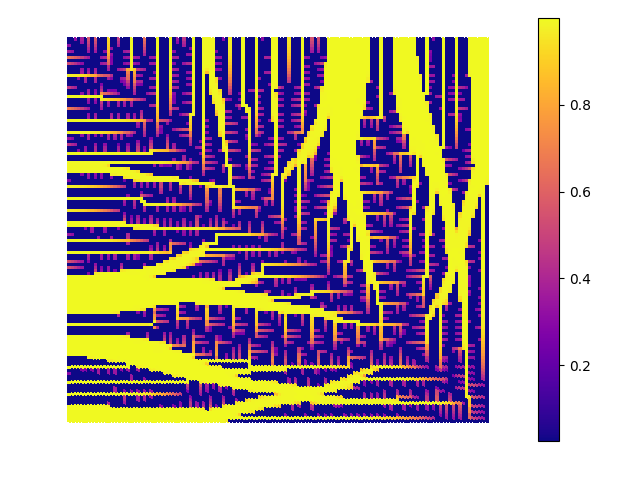}}
    \caption{Model grid refinement results, $C = 1.0$.}
    \label{fig:model_refinement_designs}
\end{figure}

\begin{figure}[]
    \centering
        \subfloat[15 iteration]{
        \label{con15}
        \includegraphics[width=0.3\columnwidth]{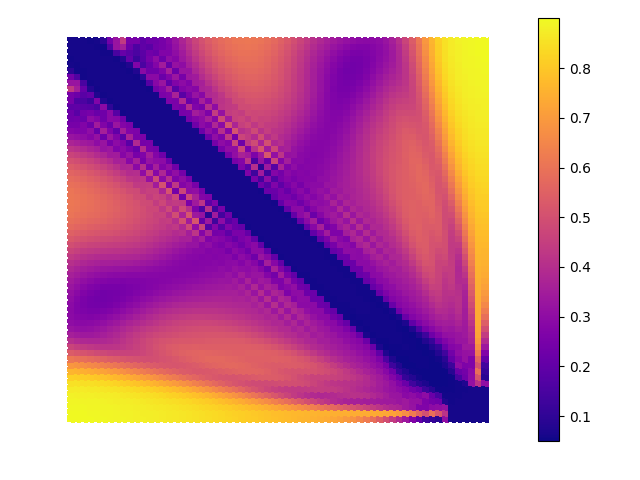}}
        \subfloat[22 iteration]{
        \label{con22}
        \includegraphics[width=0.3\columnwidth]{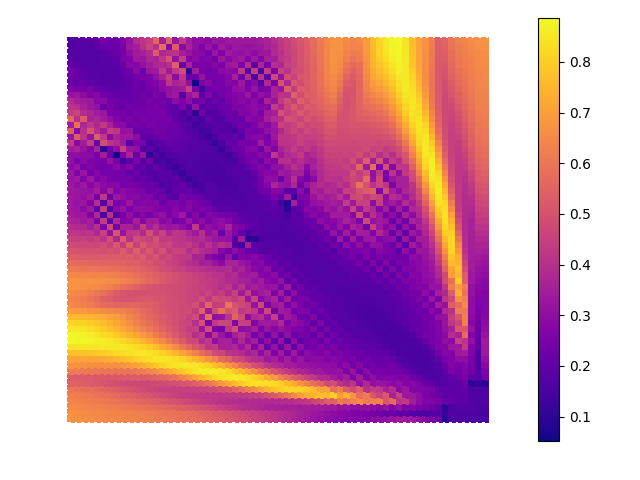}}
        \subfloat[28 iteration]{
        \label{con28}
        \includegraphics[width=0.3\columnwidth]{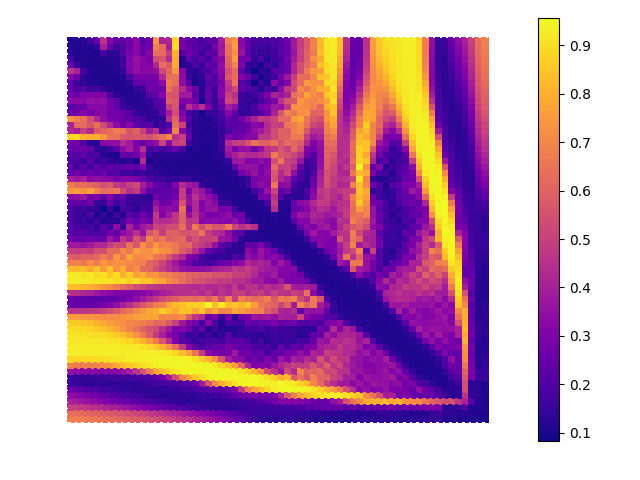}}
        
        \subfloat[47 iteration]{
        \label{con47}
        \includegraphics[width=0.3\columnwidth]{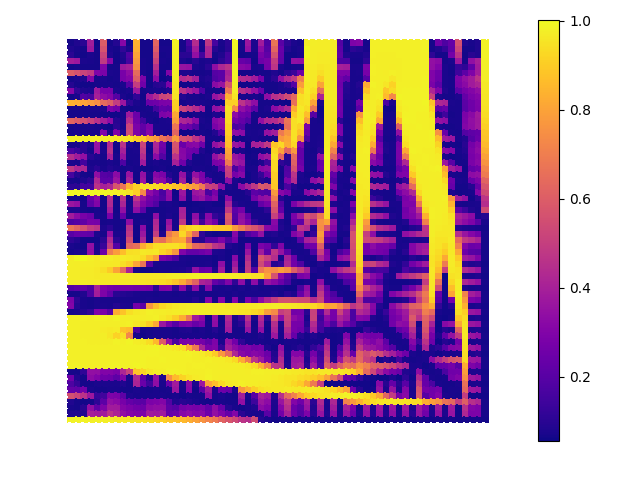}}
        \subfloat[80 iteration]{
        \label{con80}
        \includegraphics[width=0.3\columnwidth]{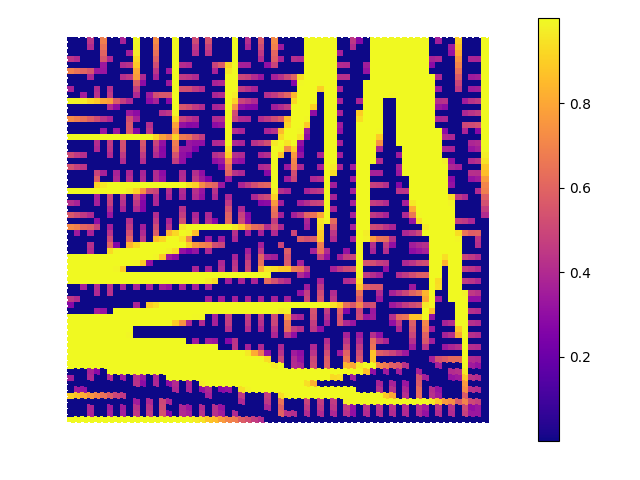}}
    \caption{Intermediate designs, $C=1.0$, $N=64$. }
    \label{fig:convergence_steps}
\end{figure}

\section*{Acknowledgement}
We deeply thank anonymous referees and Raphael Haftka for encouraging and very thorough review of our manuscript, useful citations, and very constructive criticism.

\bibliography{paper}

\end{document}